\documentclass[11pt,a4paper]{article}

\usepackage{graphicx}  
\usepackage{amsmath,amsfonts,amssymb}  
\usepackage{color}
\usepackage{url}

\newcommand{\eref}[1]{(\ref{#1})}
\newcommand{\fref}[1]{figure~\ref{#1}}
\newcommand{\Fref}[1]{Figure~\ref{#1}}
\newcommand{\frefs}[1]{figures~\ref{#1}}
\newcommand{\Frefs}[1]{Figures~\ref{#1}}
\newcommand{\sref}[1]{section~\ref{#1}}
\newcommand{\Sref}[1]{Section~\ref{#1}}

\newcommand{\eps}{\varepsilon}

\newcommand{\bbZ}{\mathbb{Z}}
\newcommand{\bbR}{\mathbb{R}}

\newcommand{\hW}{\widehat{W}}
\newcommand{\hWu}{\hW^u}

\begin{document}

\begin{center}
{\Huge Global manifold structure of a \\[2mm] continuous-time heterodimensional cycle}

\vspace*{5mm}

{\Large Andy Hammerlindl$^1$, Bernd Krauskopf$\,^2$, Gemma Mason$^2$ \\[1mm] and Hinke M.\ Osinga$^2$}
\end{center}

\vspace*{1mm}

\noindent
{$^1$ School of Mathematical Sciences, Monash University, Melbourne VIC 3800, Australia \\[1mm] 
$^2$ Department of Mathematics, The University of Auckland, Auckland 1142, New Zealand}




\vspace*{2mm}

\begin{center}
{\large June 2019}
\end{center}

\vspace*{3mm}

\begin{abstract}
A heterodimensional cycle consists of a pair of heteroclinic connections between two saddle periodic orbits with unstable manifolds of different dimensions. Recent theoretical work on chaotic dynamics beyond the uniformly hyperbolic setting has shown that heterodimensional cycles may occur robustly in diffeomorphisms of dimension at least three. We study a concrete example of a heterodimensional cycle in the continuous-time setting, specifically in a four-dimensional vector field model of intracellular calcium dynamics. By employing advanced numerical techniques, Zhang, Krauskopf and Kirk [\emph{Discr. Contin. Dynam. Syst. A\/} \textbf{32}(8) 2825--2851 (2012)] found that a heterodimensional cycle exists in this model. 

We investigate the geometric structure of the associated stable and unstable manifolds in the neighbourhood of this heterodimensional cycle, consisting of a single connecting orbit of codimension one and an entire cylinder of structurally stable connecting orbits between two saddle periodic orbits. We employ a boundary-value problem set-up to compute their stable and unstable manifolds, which we visualize in different projections of phase space and as intersection sets with a suitable three-dimensional Poincar\'{e} section. We show that, locally near the intersection set of the heterodimensional cycle, the manifolds interact as described by the theory for three-dimensional diffeomorphisms. On the other hand, their global structure is more intricate, which is due to the fact that it is not possible to find a Poincar\'{e} section that is transverse to the flow everywhere. Our results show that the abstract concept of a heterodimensional cycle arises and can be studied in continuous-time models from applications.
\end{abstract}


\section{Introduction}
\label{sec:intro}

One of the best-studied forms of chaotic dynamics is \emph{uniform hyperbolicity}, which arose from the work of Anosov and Smale in the 1960s~\cite{Anosov, Smale}. In a uniformly hyperbolic system, every orbit tends toward one of a finite number of basic sets in phase space. Each point in a basic set has well-defined directions of expansion and contraction, and tangent to these directions are stable and unstable manifolds through the point. This form of dynamics is robust, that is, any small perturbation of a uniformly hyperbolic system is still uniformly hyperbolic. Hence, it is likely that such chaotic systems may be observed in physical processes, and evidence of such behaviour has been found in a number of settings; for example, see~\cite{sturman2014, galias1997chua, rre2003explicit, viana2011fractal}. Shortly after its discovery, Smale and others conjectured that uniform hyperbolicity was the ``default'' form of chaotic dynamics. That is, they conjectured that any chaotic dynamical system could, after a small perturbation, be turned into a uniformly hyperbolic system. Consequently, most chaotic behaviour identified in a mathematical model of a physical process would be uniformly hyperbolic. Unfortunately, this hypothesis turned out to be false and other forms of robust chaotic dynamics have since been identified; see~\cite[Chapter 1]{bdv2005book}. Since these more complicated forms of non-hyperbolic chaotic dynamics also persist under perturbation, they could conceivably be observed in physical processes as well.

\begin{figure}[t!]
  \hspace*{2.4cm}
  \includegraphics{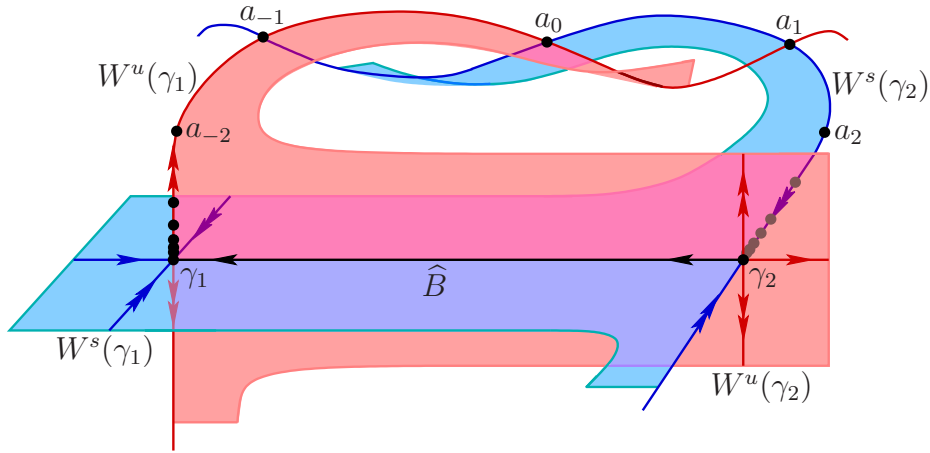}
  \caption{\label{fig:hdcycle}
    Sketch of a heterodimensional cycle in a three-dimensional discrete-time system. Here, two saddle fixed points $\gamma_1$ and $\gamma_2$ have two-dimensional manifolds $W^s(\gamma_1)$ and $W^u(\gamma_2)$ that intersect transversely in a curve $B$ and one-dimensional manifolds $W^u(\gamma_1)$ and $W^s(\gamma_2)$ that intersect in a single orbit $\left( a_k \right)_{k \in \bbZ}$. Reproduced from \cite{ZKK2012}.}
\end{figure}

One mechanism that creates non-hyperbolic chaotic dynamics is a \emph{heterodimensional cycle}~\cite{Bonatti2008, 
bdv2005book, Diaz95,yorke}. Consider a discrete-time dynamical system, a diffeomorphism, with two hyperbolic fixed (or periodic) points $\gamma_1$ and $\gamma_2$. These points form a cycle if the stable manifold $W^s(\gamma_1)$ of $\gamma_1$ intersects the unstable manifold $W^u(\gamma_2)$ of $\gamma_2$ and the stable manifold $W^s(\gamma_2)$ of $\gamma_2$ intersects the unstable manifold $W^u(\gamma_1)$ of $\gamma_1$. The cycle is heterodimensional if $\gamma_1$ has a different index from $\gamma_2$, that is, their stable manifolds have different dimensions. This requires a diffeomorphism on a phase space of dimension at least three. \Fref{fig:hdcycle} shows a sketch of a heterodimensional cycle for this case of lowest dimension, where $W^s(\gamma_1)$ and $W^u(\gamma_2)$ are two dimensional and $W^u(\gamma_1)$ and $W^s(\gamma_2)$ are one dimensional. The intersection $B$ between the surfaces $W^s(\gamma_1)$ and $W^u(\gamma_2)$ is transverse and, therefore, persists and is structurally stable under perturbation of the dynamics. The one-dimensional manifolds $W^u(\gamma_1)$ and $W^s(\gamma_2)$, on the other hand, intersect non-transversely in the heteroclinic orbit  $\left( a_k \right)_{k \in \bbZ}$. Indeed,  this intersection is of codimension one and can be destroyed by a perturbation of the system. Nevertheless, it has been proven that, near a given heterodimensional cycle as shown in \fref{fig:hdcycle}, the existence of a heterodimensional cycle is a $C^1$-open property~\cite{Bonatti2008,bdv2005book, Diaz95}. That is, the existence of any specific heterodimensional cycle is not itself a $C^1$-open property, but for any diffeomorphism $f$ with a heterodimensional cycle and any sufficiently small $\varepsilon$, there is a $C^1$-open set $U$ in the space of diffeomorphisms $\varepsilon$-close to $f$ where every element of $U$ has a heterodimensional cycle. This phenomenon is closely associated with objects called \emph{blenders}, which are invariant sets of diffeomorphisms with the property, in the present context of dimension three, that their one-dimensional stable or unstable manifold cannot be avoided by one-dimensional curves from (an open neighbourhood of) a certain direction~\cite{bcdw2016what, bd1996blender, katsutoshi, hkos-blender}. Blenders are in some sense a generalization of Smale's horseshoe construction to higher dimensions, and their robust defining property is a form of partially hyperbolic dynamics that has been used to prove the existence of heterodimensional cycles. 

Another mechanism to create non-hyperbolic chaotic dynamics is a \emph{homoclinic tangency} between a stable and an unstable manifold of a single fixed or periodic point of a diffeomorphism~\cite{pt1993hyperbolicity}. Such a non-transverse intersection may also be a robust phenomenon when the phase space of the diffeomorphism is at least three~\cite{bamon, GonchenkoMeiss2006, GonchenkoTuraev2005, GonchenkoTuraev2009, hko-wild, turaev1998}. Homoclinic tangencies and heterodimensional cycles are closely related; see~\cite[\S1.1]{cp2015essential}. A conjecture of Palis states that every dynamical system may be perturbed (in the $C^r$-topology of diffeomorphisms) to produce a system that is either uniformly hyperbolic or has a homoclinic tangency or a heterodimensional cycle~\cite{pal1991homoclinic, pt1993hyperbolicity}. Limited versions of this conjecture have been verified in a number of settings~\cite{bc2004recurrence, cp2015essential}.

The progress made in the abstract study of possible forms of chaotic dynamics has primarily focussed on diffeomorphisms, that is, on discrete-time dynamical systems. One major reason is that the phase-space dimension needed for discrete-time systems to create a particular type of chaos is always smaller than that needed for continuous-time systems. In many applications, on the other hand, mathematical models have continuous time and are given in the form of a system of ordinary differential equations (ODEs), that is, their dynamics is determined by a flow; see, for example, \cite{Edelstein, ErmentroutTerman10, GH, JacksonRadunskaya15} as entry points to the extensive literature on models from applications. Hence, it is an important question how different types of dynamics that have been found in diffeomorphisms manifest themselves in flows.

Of course, there is a well-known connection between continuous-time and discrete-time systems. Given a diffeomorphism one can construct its suspension, which is an abstract flow with an additional direction that represents the continuous time and whose time-one map is the diffeomorphism. Clearly, periodic points of the diffeomorphism correspond to periodic orbits of the flow. On the other hand, a periodic orbit of a vector field corresponds to fixed points of the associated Poincar\'{e} return map, which is obtained by introducing a (local) cross section transverse to the periodic orbit and, hence, the flow nearby. Indeed, the stability and bifurcation analysis of periodic orbits of vector fields is derived from the corresponding stability and bifurcations of the associated fixed points~\cite{Kuznetsov}. More precisely, the dynamics of a vector field locally in the neighbourhood of a periodic orbit can be described exactly in terms of the dynamics of a diffeomorphism locally in the neighbourhood of the corresponding fixed point of the associated Poincar\'{e} return map, and vice versa. However, when the dynamics is more complicated, and certainly in parameter regions where chaotic behaviour occurs, it may be necessary to consider a rather larger domain of definition of the Poincar\'{e} return map. Unless the vector field is periodically forced, that is, the flow is a suspension, it is not possible to find a section (codimension-one submanifold) to which the flow is transverse at any point. Hence, the Poincar\'{e} return map can only be defined locally in regions where the flow is transverse to the section and, moreover, the flow returns back to the section~\cite{Birkhoff17,Dullin95,Lee2008,Just00}. Especially when one is interested in chaotic dynamics associated with homoclinic and heteroclinic bifurcations of closed periodic orbits~\cite{HomburgSandstede10, rademacher05, rademacher10}, the overall picture of how their invariant manifolds intersect a codimension-one section of interest cannot be determined from the information of an equivalent diffeomorphism alone. Indeed, advanced numerical methods have an important role to play in determining the global organization of phase space, how it changes with parameters and connects with known theory~\cite{England2005, redbook}.

In this paper, we investigate a heterodimensional cycle in a four-dimensional vector field that represents a simple model of intracellular calcium oscillations. The model is based on the so-called Atri model~\cite{atri} and, in the form studied in~\cite{Zhang06}, it is given by the four differential equations
\begin{equation}
\label{eq:4DAtri}
  \left\{ \begin{array}{rcl}
       {\displaystyle \dot{c}}        &=& v, \\[2mm]
       {\displaystyle  D \, \dot{v}} &=& {\displaystyle  s \, v - \left( \alpha + k_f \, \frac{c^2}{c^2+\phi_1^2} \, n \right) 
                                  \left( \frac{\gamma \, (c_t + D \, v - s \, c)}{s} - c \right)} \\[3mm]
            & & \quad\; + k_s \, c - \eps \, (J - k_p \, c), \\[2mm]
       {\displaystyle \dot{c_t}}     &=& \eps \, (J - k_p \, c), \\[2mm]
       {\displaystyle  s \, \dot{n}} &=& {\displaystyle \frac{1}{2} \left( \frac{\phi_2}{\phi_2 + c} - n \right)}
  \end{array} \right.
\end{equation}
for the concentration of calcium $c$ inside the cell, its voltage $v$, the total calcium level $c_t$ (including that stored in the so-called endoplasmic reticulum (ER), which is a safe internal storage of this toxic substance) and the gating variable $n$ that describes transport of calcium through the membrane. The Atri model takes spatial variation into account, system~\eref{eq:4DAtri} is written in moving-frame coordinates and differentiation is with respect to the travelling-wave coordinate. Hence, periodic and solitary calcium waves are represented by periodic and homoclinic orbits of~\eref{eq:4DAtri}, respectively. We consider how the possible dynamics depends on the wave speed $s$ and the parameter $J$, which represents a given flux of calcium entering the cell from the outside in. All other parameters are constant and their values are given in Table~\ref{tab:constants}. 

\begin{table}[t!]
\begin{center}
\begin{tabular}{|c|c|c|c|c|c|c|c|c|}
  \hline
  $D$ & $\alpha$ & $k_f$ & $\phi_1$ & $\gamma$ & $k_s$ & $\eps$ & $k_p$ & $\phi_2$ \\
  \hline
  $25.0$ & $0.05$ & $20.0$ & $2.0$ & $5.0$ & $20.0$ & $0.2$ & $20.0$ & $1.0$ \\
  \hline
\end{tabular}
\caption{\label{tab:constants}
  Parameter values used for the intracellular calcium model~\eref{eq:4DAtri}.}
\end{center}
\end{table}

The Atri model~\eref{eq:4DAtri} exhibits very complicated dynamics~\cite{Zhang06} and, in particular, Zhang, Krauskopf, and Kirk~\cite{ZKK2012} realized that this system has a heterodimensional cycle between two periodic orbits with different indices. This global object was found numerically by careful computations with a new numerical method (an implementation of Lin's method) developed and presented in~\cite{ZKK2012}. To our knowledge, system~\eref{eq:4DAtri} is still the only known explicit vector field arising from an application that exhibits a heterodimensional cycle between two periodic orbits, and it has the lowest required dimension. We will refer to this heterodimensional cycle also as a PtoP cycle (for Periodic orbit to Periodic orbit); note that this short-hand notation comes from the literature on computing connecting orbits of vector fields, where heteroclinic cycles between an equilibrium and a periodic orbit (such as the one that gives rise to the chaotic attractor in the Lorenz system) are also referred to as EtoP cycles \cite{DKO15,krauskopfriess,ZKK2012}.

We present a case study of how the abstract theory of heterodimensional cycles manifests itself in a concrete vector field such as system~\eref{eq:4DAtri}. More specifically, we show how the respective global stable and unstable manifolds of the two periodic orbits intersect and give rise to the transverse and non-transverse parts of the heterodimensional cycle. We employ a state-of-the-art numerical approach based on boundary-value problem formulations~\cite{England2005, redbook} to compute and visualize selected families of orbit segments; these computations are conducted with the pseudo-arclength continuation package {\sc Auto}~\cite{AutoOrig, auto}. We show three-dimensional projections of the four-dimensional phase space to illustrate one by one how the intersection sets of the different global manifolds with a suitable codimension-one section $\Sigma$ arise. Taken together, this information makes clear how (the intersection set of) the heterodimensional cycle is generated, to what extent it represents the theory for diffeomorphisms, and what the differences are. More generally, the results presented here show that it is entirely feasible with state-of-the-art numerical methods to find and investigate objects as complex as a heterodimensional PtoP cycle in continuous-time models from applications. Hence, theoretical insights on new types of robust dynamics can now be investigated in terms of their relevance in specific contexts.

This paper is organized as follows. In \sref{sec:hetcycle} we present the heterodimensional cycle in the four-dimensional phase space of system~\eref{eq:4DAtri}. \Sref{sec:returnmap} is devoted to the manifold structure that gives rise to the PtoP cycle. The codimension-one section $\Sigma$ is introduced in \sref{sec:P}, and \sref{sec:P_PtoP} shows the intersection sets of the PtoP cycle with $\Sigma$. \Sref{sec:ws} then introduces the intersection set of the two-dimensional stable manifold and \sref{sec:wu} that of the two-dimensional unstable manifold. How the two intersection sets interact in the three-dimensional section $\Sigma$ is explained in \sref{sec:wswu}, which also contains a comparison with the case of a diffeomorphism. We end in \sref{sec:conclusions} with some conclusions and suggestions for future research.

\section{The heterodimensional PtoP cycle}
\label{sec:hetcycle}

In this paper, we explore the global manifold structure giving rise to a codimension-one heterodimensional PtoP cycle in system~\eref{eq:4DAtri}. The PtoP cycle was found in~\cite{ZKK2012} for fixed $s = s^\ast = 9.0$ at $J = J^\ast \approx 3.02661$ with an implementation of Lin's method, and its locus was then found as a curve in the $(J, s)$-plane. We start our investigation in the same way, where, as in~\cite{ZKK2012}, all other parameters are as given in Table~\ref{tab:constants}.

\begin{figure}[t!]
  \hspace*{1.5cm}
  \includegraphics{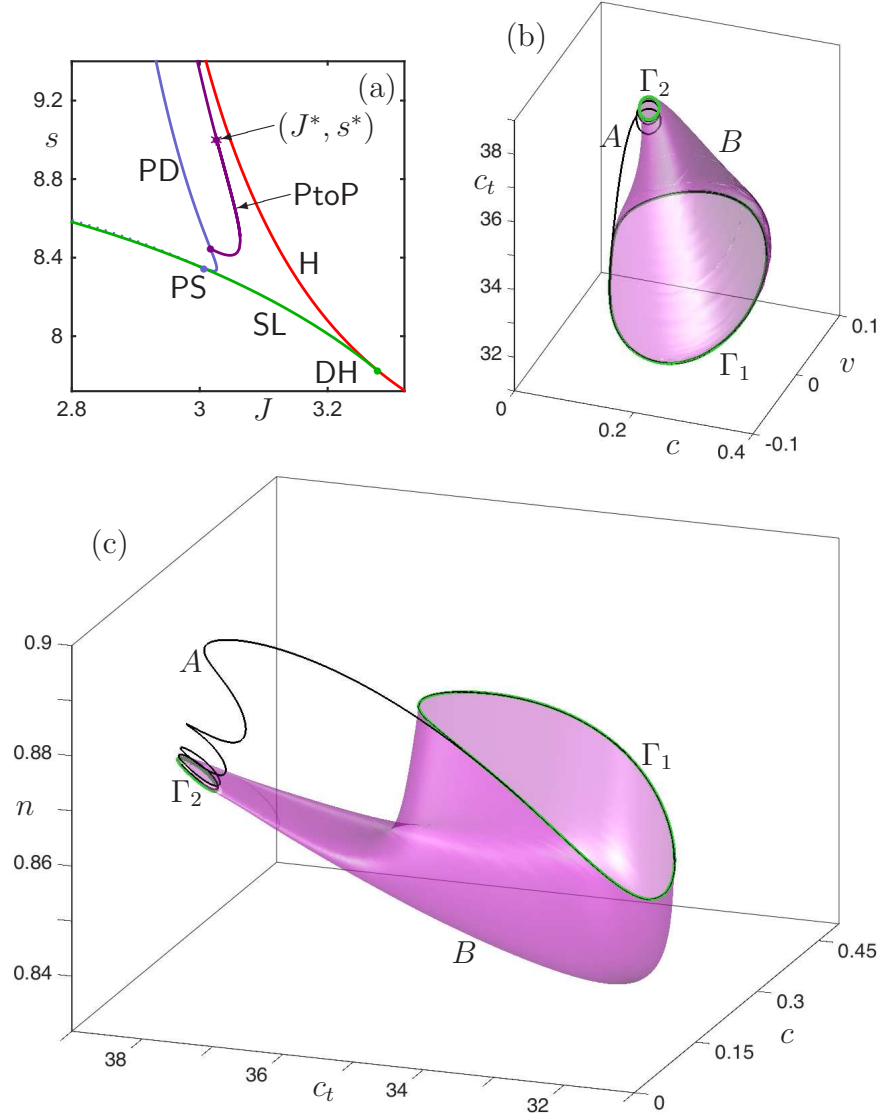} 
  \caption{\label{fig:hetcycle}
    Panel (a) shows the locus {\sf PtoP} (purple curve) of the heterodimensional cycle of system~\eref{eq:4DAtri} in the $(J, s)$-plane, relative to the loci of Hopf bifurcation {\sf H} (red curve), of saddle-node bifurcation of limit cycles {\sf SL} (green curve) ending on {\sf H} at the point {\sf DH}, and of period-doubling bifurcation {\sf PD} (blue curve); the curve {\sf PD} is tangent to {\sf SL} at the point {\sf PS}, to the left of which {\sf PD} is dotted. The heterodimensional PtoP cycle for the indicated point $(J^\ast, s^\ast) = (3.02661, 9.0)$ on {\sf PtoP} is shown in projections into $(c, v, c_t)$-space in panel~(b) and into $(c, c_t, n)$-space in panel~(c). It consists of a unique (and non-transverse) connecting orbit $A$ (black curve) from $\Gamma_1$ to $\Gamma_2$ (green curves) and a two-dimensional topological cylinder $B$ (purple surface) of trajectories from $\Gamma_2$ to $\Gamma_1$.}
\end{figure}

\Fref{fig:hetcycle} shows the respective part of the $(J, s)$-plane in panel~(a) and two three-dimensional projections of the PtoP cycle at $(J^\ast, s^\ast)$ in panels~(b) and~(c). In \fref{fig:hetcycle}(a), two periodic orbits $\Gamma_1$ and $\Gamma_2$ bifurcate from the curve {\sf SL} of saddle-node bifurcation of periodic orbits. These two periodic orbits are both hyperbolic, of saddle type and of different index in the region that is bounded by {\sf SL} and the curves {\sf H} of Hopf bifurcation, where $\Gamma_2$ bifurcates, and {\sf PD} of period-doubling bifurcation, where $\Gamma_1$ bifurcates. Notice that the curve {\sf SL} ends on {\sf H} at a degenerate Hopf point {\sf DH} of codimension two, while the curve {\sf PD} is tangent to {\sf SL} at a codimension-two saddle-node-period-doubling point {\sf PS}. The heterodimensional cycle exists along the curve {\sf PtoP}, which ends at a point (not labelled) on the curve {\sf PD}. 

The point $(J^\ast, s^\ast)$ is marked, and the corresponding heterodimensional PtoP cycle is shown in
panels (b) and~(c) of \fref{fig:hetcycle}. More specifically, these images show the PtoP cycle between the two periodic orbits $\Gamma_1$ and $\Gamma_2$, which have Floquet multipliers
\begin{displaymath}
  \lambda^u_1 \approx 9.6 \times 10^4, \quad 
  \lambda^s_1 \approx 2.8 \times 10^{-1}, \quad \mbox{and} \quad 
  \lambda^{ss}_1 \approx 7.7 \times 10^{-3}, 
\end{displaymath} 
and
\begin{displaymath}
  \lambda^{uu}_2 \approx 7.0 \times 10^3 \quad 
  \lambda^u_2 \approx 1.3, \quad \mbox{and} \quad 
  \lambda^s_2 \approx 3.4 \times 10^{-1}, 
\end{displaymath} 
respectively. Hence, $\Gamma_1$ has a two-dimensional unstable manifold $W^u(\Gamma_1)$ and a three-dimensional stable manifold $W^s(\Gamma_1)$, while $\Gamma_2$ has a three-dimensional unstable manifold $W^u(\Gamma_2)$ and a two-dimensional stable manifold $W^s(\Gamma_2)$. The three-dimensional manifolds $W^s(\Gamma_1)$ and $W^u(\Gamma_2)$ intersect transversely, as expected. The two-dimensional manifolds $W^s(\Gamma_2)$ and $W^u(\Gamma_2)$, on the other hand, have a non-transverse intersection at the special parameter pair $(J^\ast,\, s^\ast)$. We denote these non-empty intersections by
\begin{displaymath}
  A : = W^s(\Gamma_2) \cap W^u(\Gamma_1), \quad \mbox{and} \quad B := W^s(\Gamma_1) \cap W^u(\Gamma_2).
\end{displaymath}
The intersection set $A$ is one dimensional and consists of a single trajectory from $\Gamma_1$ to $\Gamma_2$. The intersection set $B$ is two dimensional and forms a topological cylinder consisting of a one-parameter family of trajectories from $\Gamma_2$ to $\Gamma_1$. Note the strong contraction toward $\Gamma_1$ in backward time and the much slower approach toward $\Gamma_2$ in forward time, which is particularly visible in panel~(c); this difference in speed is due to the fact that the ratio $\lambda^s_2 / \lambda^{ss}_1$ is of the order $10^2$. The curve  {\sf PtoP} in the $(J, s)$-plane in panel~(a) was computed by parameter continuation of the non-transverse, codimension-one connecting orbit $A$, while checking that the cylinder $B$, that is, the entire heterodimensional cycle, also exists along {\sf PtoP}.

\subsection{Numerical approach to finding the PtoP cycle}
\label{sec:computation}

As mentioned in the introduction, the computation of the heterodimensional cycle for system~\eref{eq:4DAtri} requires the advanced numerical approach from~\cite{ZKK2012}. Its major element is that it is based on Lin's method~\cite{Lin}, which has the important property that it defines a smooth test function --- thereby  guaranteeing (under very mild starting conditions) that a particular homoclinic or heteroclinic orbit is found if it exists. We refer to~\cite{ZKK2012} for the full details, but give a brief overview here of the steps involved. 

Numerical methods for approximating connecting orbits are not new; see~\cite{beyn90} and references therein. Initially, the interest was in computing connections between equilibria. Depending on the dimensions of the respective eigenspaces, such connections may exist only for specific parameter values, which gives an added computational challenge. The connecting orbit is represented by an orbit segment truncated at both ends; it is assumed to start close to an equilibrium in its unstable eigenspace and end close to the same or other equilibrium in its stable eigenspace. The same approach can also be used for connections involving periodic orbits, where the begin or end point of the orbit segment is required to lie in the linear eigenbundle associated with the stable or unstable Floquet multipliers. The connecting orbit is then found as the orbit segment that solves a two-point boundary value problem (2PBVP)~\cite{beyn90, beyn94}; this is typically done in a continuation setting, where the specific connecting orbit is detected as a special solution from a family with a known solution. 

To achieve a well-posed and numerically robust set-up, Beyn~\cite{beyn94} popularized the use of projection boundary conditions for the begin and end points of the orbit segment, which provides full error control. Various implementations exist, involving the continuation of invariant subspaces~\cite{dieci04} or the (adjoint) variational equation~\cite{vanVoorn08, vanVoorn09}. The examples in these papers include computation of so-called EtoP and PtoP connections, that is, heteroclinic orbits between an equilibrium (E) and a periodic orbit (P), and homoclinic or heteroclinic orbits between one or two periodic orbits (PtoP), respectively; note that the example of a PtoP connection from~\cite{dieci04} is between two periodic orbits with the same number of stable Floquet multipliers, so they are not part of a heterodimensional cycle. The continuation of the 2PBVP is started from a solution found by shooting~\cite{dieci04} or via a homotopy step~\cite{vanVoorn08, vanVoorn09}; the latter is also the standard approach for the software package {\sc HomCont}~\cite{homcont} that computes homoclinic and heteroclinic orbits of equilibria. 

The set-up based on Lin's method provides a natural way to find a first connecting orbit segment. The implementation for connections between equilibria was presented in~\cite{Oldeman03} and the approach was extended to connections involving periodic orbits in~\cite{krauskopfriess,ZKK2012}. The idea is to split the connecting orbit into two halves and require that both halves have an end point in a suitably chosen codimension-one section $\Pi$ that is (locally) transverse to the flow. At the same time, it is required that the two orbit segments satisfy suitable projection boundary conditions at their other end points. Typically, the two end points in $\Pi$ will not be the same, but if they are then a connecting orbit has been found. The key insight of Lin's method is that the gap between the two end points in $\Pi$, the so-called Lin gap, can be restricted to lie in a particular subspace $Z \subset \Pi$ of a dimension related to the codimension of the connecting orbit. In other words, the Lin gap can be closed in a component-wise manner with respect to a chosen basis for $Z$ while adjusting the start or end point of one of the orbit segments within the required linear eigenspace. The computational set-up requires continuation of two orbit segments (double the system dimension), but this is outweighed by the major advantage that the two-segment 2PBVP is well posed even when the parameters are such that the particular connecting orbit does not exist; as long as $\Pi$ is transverse to the flow for the family of orbit segments under consideration, a first solution to the 2PBVP can be found relatively easily and a possible special orbit segment with closed Lin gap can be searched for in the family. 

The codimension-one connecting PtoP orbit  $A = W^u(\Gamma_1) \cap  W^s(\Gamma_2)$ from $\Gamma_1$ to $\Gamma_2$ in panels~(b) and~(c) of \fref{fig:hetcycle} has been found in this way; here the Lin gap is of dimension one and closing it requires changing a system parameter, which was $J$ in this case. The corresponding pair of connecting orbit segments can then be continued in system parameters with the Lin gap remaining closed, or as a single orbit segment after concatenation of the pair into a single orbit segment. The curve {\sf PtoP} in \fref{fig:hetcycle}(a) was computed in this way.

The cylinder $B = W^u(\Gamma_2) \cap  W^s(\Gamma_1)$ from $\Gamma_2$ to $\Gamma_1$ in panels~(b) and~(c) of \fref{fig:hetcycle} can be computed with effectively the same approach. One first finds a pair of orbits with projection boundary conditions in $W^s(\Gamma_1)$ and $W^u(\Gamma_2)$, respectively, with their other end points in $\Pi$. Again, a one-dimensional Lin gap can then be defined in $\Pi$, but to close it one does not need to change system parameters in this case. Rather, varying the distance from one of the periodic orbits closes the gap. Once a first connecting orbit has been found, the cyclinder can be swept out as a one-parameter familly while varying the distance from one of the periodic orbits (over one fundamental domain of the local return); see~\cite{ZKK2012} for the more details.

\section{Global manifold structure of the PtoP cycle}
\label{sec:returnmap}

Our goal is to investigate how the heterodimensional PtoP cycle $\Gamma_1 \cup A \cup \Gamma_2 \cup B$ in the four-dimensional system~\eref{eq:4DAtri} arises from the intersections of the two-dimensional manifolds $W^u(\Gamma_1)$ and $W^s(\Gamma_2)$ and the three-dimensional manifolds of $W^u(\Gamma_2)$ and $W^s(\Gamma_1)$, respectively. Moreover, we wish to make the connection between this heterodimensional PtoP cycle in a flow and the abstract theory for diffeomorphisms. The challenge is to compute and visualize parts of interest of invariant manilofds that provide the insights that we seek. To this end, we again find and continue orbit segments that are defined by suitable 2PBVPs. Indeed, this approach is very versatile and accurate when it comes to computing two-dimensional (sub)manifolds. As before, the key idea is to impose a projection boundary condition near one of the periodic orbits, while there are several choices for the boundary condition at the other end point; we refer to~\cite{redbook} for more details and examples, and to~\cite{surveymans} for an overview of different methods. In the present context, such computations can be started from connecting orbits, for example, by relaxing the requirement of a connection and sweeping out a relevant part of the respective invariant manifold up to a chosen section.

\subsection{The three-dimensional Poincar\'{e} section}
\label{sec:P}


We now consider the intersection sets of relevant invariant objects with a suitably chosen section. Recall that, in general, it is not possible to find a so-called global Poincar\'e section --- to which the flow is  transverse everywhere and to which all points return under the flow. Unless the system is periodically forced, any three-dimensional section $\Sigma$ will have a codimension-one tangency locus $C \subset \Sigma$ along which the flow, which we refer to as $\varphi^t$, is tangent to $\Sigma$.  As a result, a return map $f: \Sigma \to \Sigma$ can only be defined locally on a subset of $\Sigma \setminus C$ with the additional condition that the images under $\varphi^t$ also exist  \cite{Dullin95,Lee2008,Just00}. The maximal domain $U$ of definition of a return map given by $\varphi^t$ is then 
\begin{displaymath}
  U := \{x  \in  \Sigma \setminus C \;|\; \exists \, t > 0 \mbox{ such that } \varphi^t(x)  \in  \Sigma \setminus C \}.
\end{displaymath}
Note that both $U$ and $f(U)$ are open subsets of $\Sigma$. Any point on the boundary of $U$ either lies on the tangency locus $C$ or, as it is approached from within $U$, the return time to $\Sigma$ tends to infinity. Note that the flow direction on either side of $C$ typically changes sign with respect to the normal to $\Sigma$; hence, one can split $U$ into two disjoint subsets of $\Sigma$, denoted $U^+$ and $U^-$, which are separated by $C$. It is natural to define the return maps $f^+$ and $f^-$ on these subsets $U^+$ and $U^-$, respectively; note that they are generically given by the second return to $\Sigma$, that is, by the second iterate $f^2$ of $f$. In particular, the intersection points of a periodic orbit with $\Sigma$ will then be fixed points of $f^+$ and $f^-$, which both are then the usual local Poincar\'e maps near the periodic orbit.

Indeed, the Atri system~\eref{eq:4DAtri} is not periodically forced, and so any three-dimensional section will have a tangency locus. However, the system has the nice property that we can choose a hyperplane that is transverse to the heterodimensional cycle, that is, it is transverse to $\Gamma_1$ and $\Gamma_2$ as well as to $A$ and $B$. The geometric idea is to choose a section that contains a joint curve of rotation of both $\Gamma_1$ and $\Gamma_2$ and, throughout, we consider here the three-dimensional section 
\begin{displaymath}
  \Sigma := \{ (c,v, c_t,n) \in \bbR^4\ | \ c = 0.15 \}.
\end{displaymath}
(Note that this is not the section $\Pi$ used for finding $A$ and $B$, which was chosen such that it intersects neither $\Gamma_1$ nor $\Gamma_2$ and separates these two periodic orbits.) Moreover, equation $\dot{c} = v$ of~\eref{eq:4DAtri} implies that the tangency locus of $\Sigma$ is the two-dimensional plane 
\begin{displaymath}
  C := \{ (c,v, c_t,n) \in \bbR^4\ | \ c = 0.15 ,\, v = 0\},
\end{displaymath}
which divides $\Sigma$ into two open sets, $\Sigma^+$ where $v > 0$ and  $\Sigma^-$ where $v < 0$.

\subsection{The PtoP orbit in the section $\Sigma$}
\label{sec:P_PtoP}

\begin{figure}[t!]
  \hspace*{1.5cm}
  \includegraphics{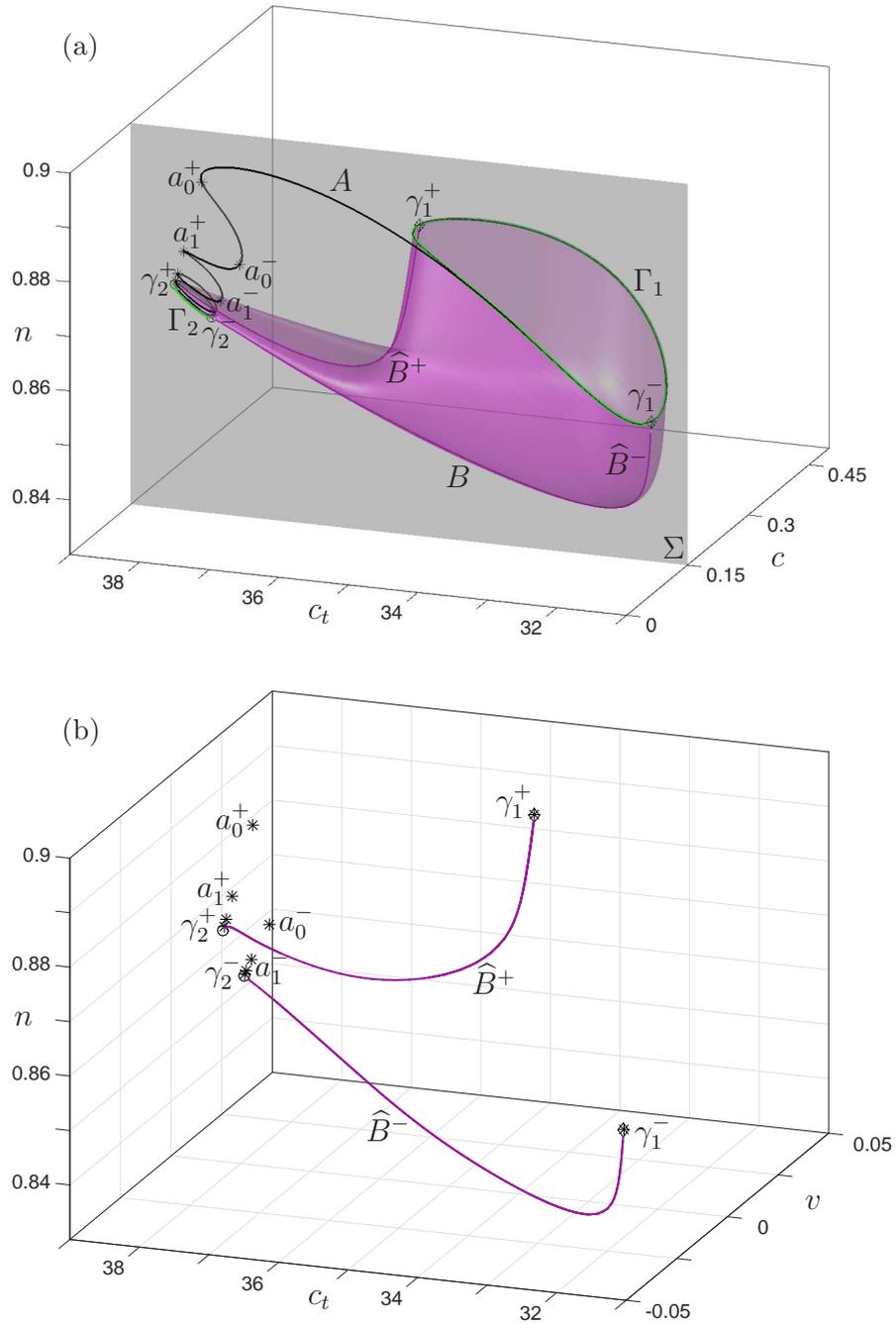} 
  \caption{\label{fig:Pcycle}
    The heterodimensional PtoP cycle, shown in panel~(a) in projection onto $(c, c_t, n)$-space with the section $\Sigma$ (grey plane) defined by $c = 0.15$, while panel~(b) shows its intersection sets in $\Sigma$. The periodic orbits $\Gamma_1$ and $\Gamma_2$ (green curves) intersect $\Sigma$ in the points $\gamma_1^\pm$ and $\gamma_2^\pm$, respectively; the connecting orbit $A$ (black curve) intersects $\Sigma$ in points $a^\pm_k$ marked by $\ast$ (these are extremely close to $\gamma^\pm_1$ for $k \leq -1$); and the cylinder $B$ (purple surface) intersects $\Sigma$ in two (purple) curves $\widehat{B}^\pm$.}
\end{figure}

\Fref{fig:Pcycle}(a) shows the heterodimensional PtoP cycle with the section $\Sigma$ and the corresponding intersection sets in projection onto $(c, c_t, n)$-space; note that $\Sigma$ appears two dimensional in this projection and compare with \fref{fig:hetcycle}(c). The intersection sets of the PtoP cycle are shown in \fref{fig:Pcycle}(b) in the three-dimensional section $\Sigma$, that is, in $(v, c_t,n)$-space with $c = 0.15$. The periodic orbits $\Gamma_1$ and $\Gamma_2$ intersect $\Sigma^+$ and $\Sigma^-$ in the points $\gamma_1^\pm$ and $\gamma_2^\pm$, respectively. The cylinder $B$ intersects $\Sigma$ in two curves $\widehat{B}^+ \subset \Sigma^+$ from $\gamma_2^+$ to $\gamma_1^+$ and $\widehat{B}^- \subset \Sigma^-$ from $\gamma_2^-$ to $\gamma_1^-$. Finally, the intersection set of the connecting orbit $A$ consists of two sequences
\begin{displaymath}
  \left( a_k^- \right)_{k \in \bbZ} := A \cap \Sigma^- \quad \mbox{and} \quad
   \left( a_k^+ \right)_{k \in \bbZ} := A \cap \Sigma^+,
\end{displaymath}
with the property that $f^\pm(a_k^\pm) = a_{k+1}^\pm$. Notice that, owing to the strong contraction toward $\Gamma_1$ in backward time, all $a_{k}^\pm$ for $k < 0$ are extremely close to $\gamma_1^\pm$, respectively.

\Fref{fig:Pcycle} shows that $\Sigma$ is indeed transverse to the entire heterodimensional cycle. However, this does not imply that all trajectories that make up the manifolds $W^u(\Gamma_1)$, $W^s(\Gamma_1)$ $W^u(\Gamma_2)$ and $W^s(\Gamma_2)$ are transverse to $\Sigma$. As a consequence, and as we will show, the intersection sets of these manifolds with $\Sigma$ can cross the tangency locus $C$ and may consist of infinitely many disjoint branches. More specifically, since $\Gamma_1$ and $\Gamma_2$ are transverse to $\Sigma$, the local return maps $f^\pm$ are diffeomorphisms near $\gamma_1^\pm$ and $\gamma_2^\pm$ and the points are hyperbolic fixed points under $f^\pm$, respectively. Hence, the Stable Manifold Theorem~\cite{PdM} guarantees the existence and smoothness of unique stable and unstable manifolds for these fixed points. However, their existence as unique smooth manifolds is guaranteed only locally. If we consider them globally, the stable and unstable manifolds of $\gamma_1^\pm$ and $\gamma_2^\pm$ --- that is, the intersection sets of the respective manifolds of $\Gamma_1$ and $\Gamma_2$ with $\Sigma$ --- may cross the tangency locus $C$ and/or consist of disjoint branches, much as invariant manifolds of endomorphisms~\cite{England2005, eko-ijbc2005, ho-ijbc2008, kop-siads2007}.

In spite of this fundamental difficulty, we now proceed with computing and visualising how the intersection sets $\left( a_k^\pm \right)_{k \in \bbZ}$ of $A$ and $\widehat{B}^\pm$ of $B$ arise from the interactions between the respective stable and unstable manifolds of the periodic orbits $\Gamma_1$ and $\Gamma_2$. To this end, we compute and discuss one by one (relevant parts of) two-dimensional manifolds in the four-dimensional phase space of system~\eref{eq:4DAtri} that lead to one-dimensional intersection curves in the three-dimensional section $\Sigma$. The final overall picture in the section $\Sigma$ that we obtain in this way allows us to relate the geometry of objects of the flow to the existing theory on heterodimensional cycles of diffeomorphisms \cite{Bonatti2008,katsutoshi}, as visualized in \fref{fig:hdcycle}.

\subsection{The intersection sets $\hW^{s, \pm}(\Gamma_2)$ and $\hW^{ss, \pm}(\Gamma_1)$}
\label{sec:ws}

By definition of the heterodimensional connection $A$, the points $\left( a_k^\pm \right)_{k \in \bbZ} = A \cap \Sigma$ must lie on the intersection set $\hW^{s,\pm}(\Gamma_2)$ of $W^s(\Gamma_2)$ with $\Sigma^\pm$. We can utilize this property to find the curves of $\hW^{s,\pm}(\Gamma_2)$ locally near points of the set $\left( a_k^\pm \right)_{k \in \bbZ}$. For example, the point $a^+_0 \in \Sigma^+$ corresponds to the part of the orbit segment representing $A$ that starts at $a^+_0$ and ends very close to $\Gamma_2$ in the linear approximation of $W^s(\Gamma_2)$. This part of the orbit segment can be selected and then continued by considering it as a solution of the 2PBVP with one end point in $\Sigma^+$ and the other still satisfying the projection boundary condition near $\Gamma_2$. 

\begin{figure}[t!]
  \hspace*{1.5cm}
  \includegraphics{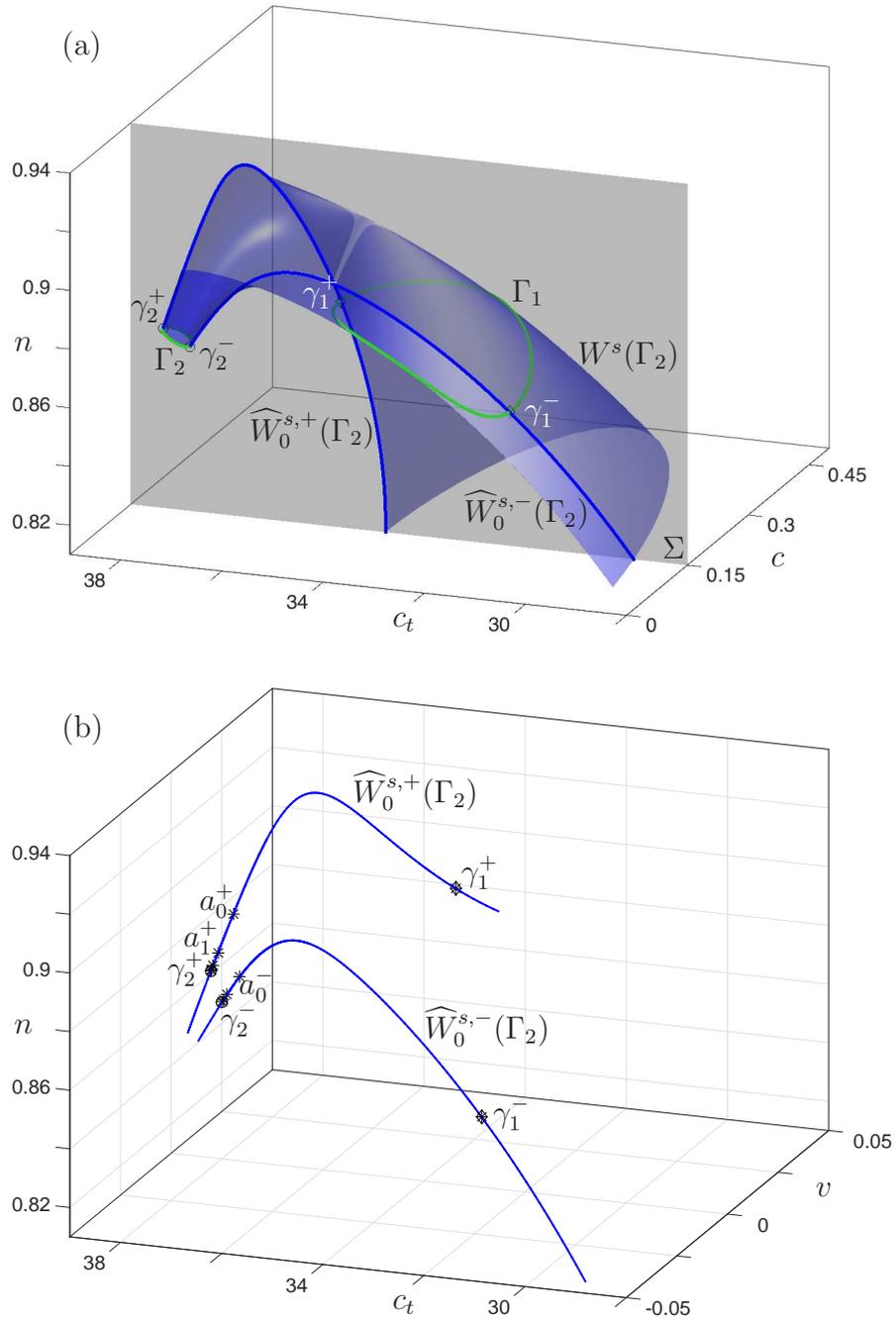} 
  \caption{\label{fig:Ws}
    The two-dimensional stable manifold $W^s(\Gamma_2)$ (blue surface) intersects $\Sigma$ (grey plane) in the two primary intersection curves $\hW_0^{s,\pm}(\Gamma_2)$ (blue curves). Panel~(a) shows, in projection onto $(c, c_t, n)$-space, the section $\Sigma$ and the side of $W^{s}(\Gamma_2)$ that comes very close to $\Gamma_1$ (green curve); panel~(b) shows in $\Sigma$ the intersection sets $\hW_0^{s,\pm}(\Gamma_2)$, $\gamma^\pm_1$,  $\gamma^\pm_1$ and $a^\pm_k$.}
\end{figure}

\Fref{fig:Ws} shows that the continuation of this orbit segment starting from $a^+_0$ gives a single curve in $\Sigma^+$ that connects all the points $ a_k^+$ for $k > -1$ to $\gamma_2^+$. Hence, this curve is (one side of) the local stable manifold of $\gamma_2^+$ under $f^+$, which we also refer to as the \emph{primary curve} $\hW_0^{s,+}(\Gamma_2)$ in $\hW^{s,+}(\Gamma_2) \subset \Sigma^+$. Panel~(a) shows, in projection onto $(c, c_t, n)$-space, how $\hW_0^{s,+}(\Gamma_2)$ arises as the intersection of $W_0^{s,+}(\Gamma_2)$ with $\Sigma$; shown are $\Gamma_1$, $\Gamma_2$, $\Sigma$ and the side of $W^s(\Gamma_2)$ that contains the connecting orbit $A$ (which is not shown). Panel~(b) shows how the primary curve $\hW_0^{s,+}(\Gamma_2)$ starts at $\gamma_2^+$, goes through the points $a_k^+$ for $k > -1$ and then passes $\gamma_1^+$ very closely, namely through $a_{-1}^+$; the latter point is the asterisk (not labelled) that lies practically on the diamond representing $\gamma_1^+$; this close passage is explained by the fact that $\Gamma_1$ has the very large unstable Floquet multiplier $\lambda^u_1 \approx 9.6 \times 10^4$. 

Similarly, starting from $a^-_0$, we find the primary curve $\hW_0^{s,-}(\Gamma_2) \subset \hW^{s,-}(\Gamma_2)$ that contains $\gamma_2^-$ as well as $a_k^-$ for $k > -1$ and is, hence, the local stable manifold of $\gamma_2^-$ under $f^-$. As \fref{fig:Ws} shows, both primary curves arise from the fact that the surface $W_0^{s,-}(\Gamma_2)$ is a cylinder transverse to $\Sigma$. Both curves are entirely contained in their respective half-spaces $\Sigma^\pm \subset \Sigma$ and appear to extend to infinity after passing through $a^\pm_{-1}$, respectively. We remark that starting the continuation from any of the points $a^\pm_k$ with $k \geq -1$ results in the computation of the same two primary curves $\hW_0^{s,\pm}(\Gamma_2)$; notice that we also show, in panel~(b), short first pieces of the other side of the primary curves $\hW_0^{s,\pm}(\Gamma_2)$, which also appear to extend to infinity. However, neither side of these curves contain any of the points $a^\pm_k$ for $k \leq -2$.

The $\lambda$-lemma~\cite{PdM} implies that $W^s(\Gamma_2)$ accumulates on the two-dimensional strong stable manifold $W^{ss}(\Gamma_1)$ of $\Gamma_1$. Since the two primary curves $\hW_0^{s,\pm}(\Gamma_2)$ extend to infinity after passing through the points $a_{-1}^\pm$, there must be other curves in the intersection set $\hW^{s,\pm}(\Gamma_2)$ that go through the points $a_k^\pm$ for $k \leq -2$; these curves must accumulate near the points $\gamma_1^\pm$ on the two primary curves  $\hW_0^{ss,\pm}(\Gamma_1)$ through the points $\gamma_1^\pm$. Additional curves in $\hW^{s,\pm}(\Gamma_2)$ can again be found by continuation of orbit segments in $W^s(\Gamma_2)$ that now start at the points $a_k^\pm$ for $k \leq -2$; we denote these successive curves by $\hW_{1}^{s,\pm}(\Gamma_2)$, $\hW_{2}^{s,\pm}(\Gamma_2)$ and so on.

\begin{figure}[t!]
  \hspace*{1.5cm}
  \includegraphics{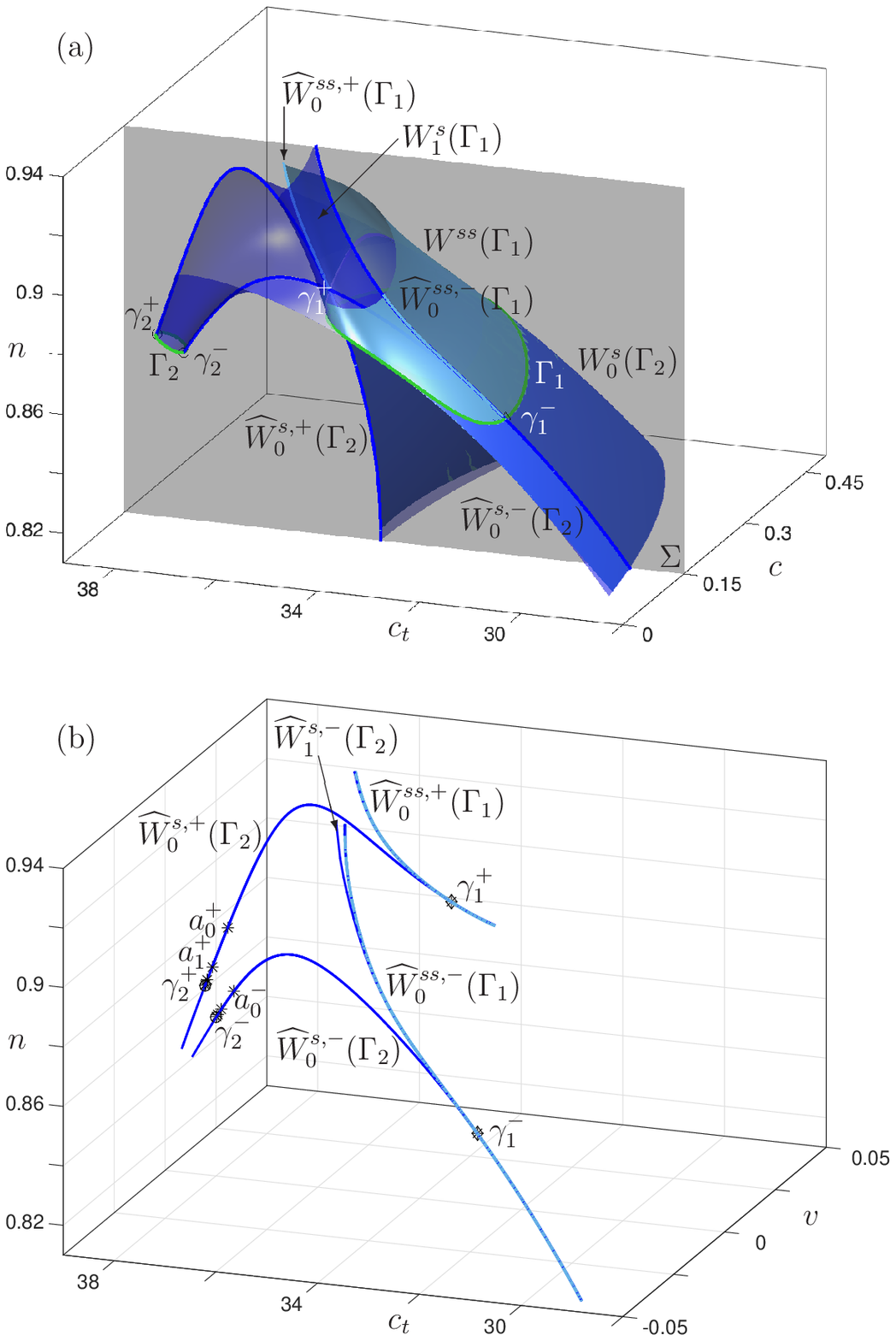} 
  \caption{\label{fig:Wss}
   The stable manifold $W^s(\Gamma_2)$ (blue) returns to $\Sigma$ (grey plane) in backward time creating additional intersection curves, of which the first, $\hW_{-1}^{s,\pm}$, is shown. These backward-time returns accumulate very fast onto the intersection set $\hW^{ss, \pm}$ (cyan curve) with $\Sigma$ of the two-dimensional strong stable manifold $W^{ss}(\Gamma_1)$ (cyan surface). Panel~(a) shows a projection onto $(c, c_t, n)$-space, and panel~(b) shows the respective intersection sets in $\Sigma$; compare with \fref{fig:Ws}.}
\end{figure}

\Fref{fig:Wss} illustrates the accumulation of $W^s(\Gamma_2)$ onto $W^{ss}(\Gamma_1)$. Panel~(a) shows, in projection onto $(c, c_t, n)$-space, $\Gamma_1$, $\Gamma_2$ and $\Sigma$, the two first pieces $W_0^s(\Gamma_2)$ and $W_1^{s}(\Gamma_2)$ of $W^{s}(\Gamma_s)$, and the part of $W^{ss}(\Gamma_1)$ that gives the intersection curve $\hW_{1}^{s,+}(\Gamma_2)$. The surface $W_1^{s}(\Gamma_2)$ is very close to $W^{ss}(\Gamma_1)$, which is illustrated in panel~(b) by the corresponding intersection curves in $\Sigma$. Due to the strong contraction in backward time, only the curve $\hW_{1}^{s,-}(\Gamma_2)$ through $a^-_{-2}$ can be seen. Already the curve $\hW_{2}^{s,-}(\Gamma_2)$ is so close to $\hW_{0}^{ss,-}(\Gamma_1)$ that they are indistinguishable in \fref{fig:Wss}(b); similarly, $\hW_{1}^{s,+}(\Gamma_2)$ through $a^+_{-2}$ and $\hW_{0}^{ss,+}(\Gamma_1)$ are indistinguishable on the scale of this figure.

We computed the curves $\hW_{j}^{s,+}(\Gamma_2)$ and $\hW_{j}^{s,-}(\Gamma_2)$ for $j \geq 1$ to very long arclengths, and found that they lie on either side of the tangency locus and extend to infinity on both sides of the respective points in $a_k^\pm$ for $k \leq -2$. Therefore, we conclude that these infinitely many curves in $\hW^{s,+}(\Gamma_2)$, which map to each other under the inverse of $f^\pm$, do not connect to form single curve with $\hW_{0}^{s,\pm}(\Gamma_2)$.

\subsection{The  intersection sets $\hW^{u}(\Gamma_1)$ and $\hW^{uu}(\Gamma_2)$}
\label{sec:wu}

By definition of the heterodimensional cycle $A$, the points in $\left( a_k^\pm \right)_{k \in \bbZ} \in A \cap \Sigma$ must also lie on the intersection of $W^u(\Gamma_1)$ with $\Sigma$. Analogously to the computations in the previous section, we now compute the intersection curves in the intersection set $\hW^{u}(\Gamma_1)$ by continuation of orbit segments that start very close to $\Gamma_1$, in the linear approximation of $W^u(\Gamma_1)$, and end in $\Sigma$, starting initially at one of the points in $\left( a_k^\pm \right)_{k \in \bbZ} \subset \Sigma$. 

\begin{figure}[t!]
  \hspace*{1.5cm}
  \includegraphics{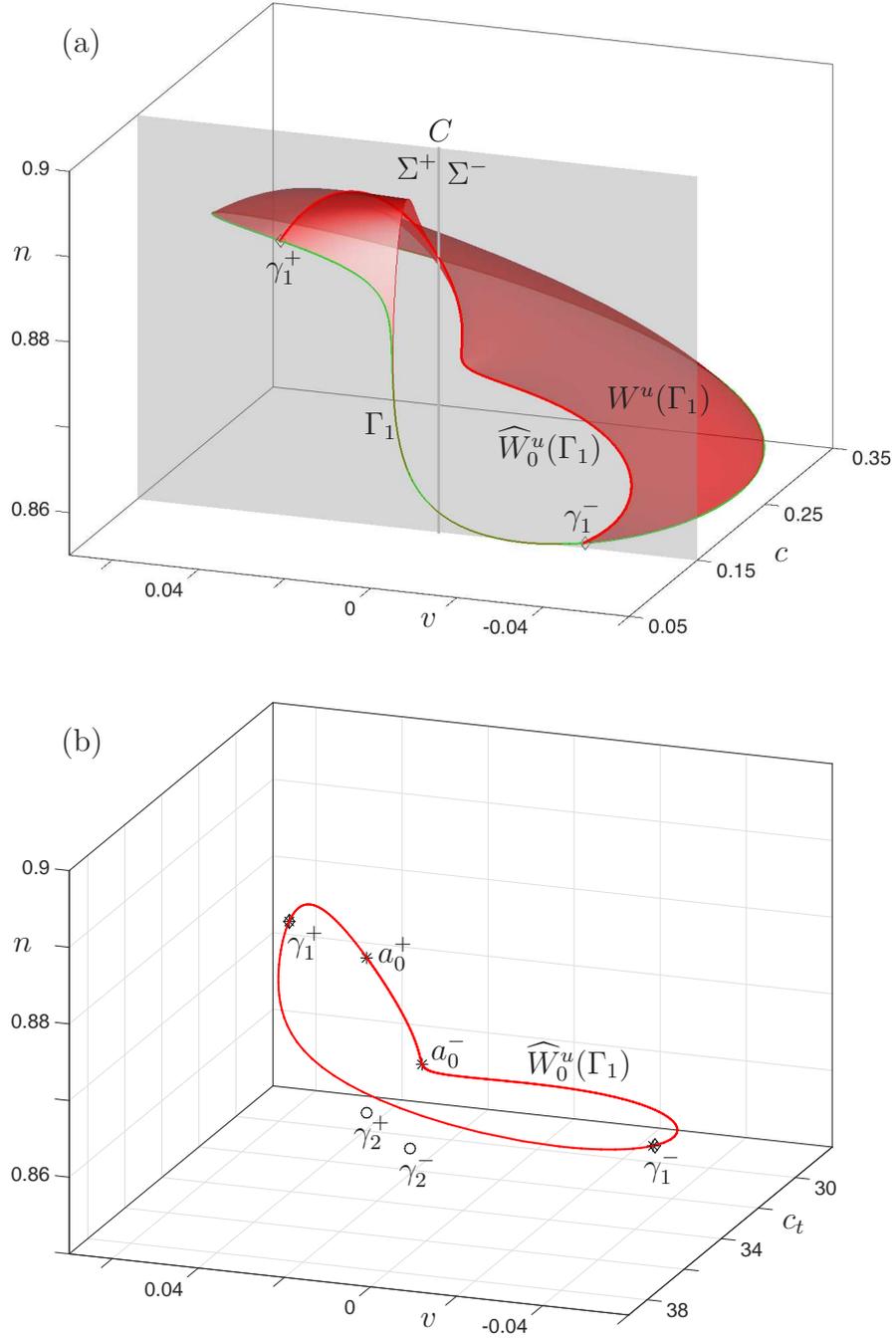} 
  \caption{\label{fig:Wu0}
    The two-dimensional unstable manifold $W^{u}(\Gamma_1)$ (red surface) intersects  $\Sigma$ (grey plane) in the primary curve $\hWu_0(\Gamma_1)$ that contains the two points $\gamma_1^\pm$ and crosses the tangency locus $C$ in $\Sigma$ twice. Panel~(a) shows, in projection onto $(c, v, n)$-space, the part of $W^{u}(\Gamma_1)$ between $\Gamma_1$ (green curve) and the arc of $\hWu_0(\Gamma_1)$ (red curve) in $\Sigma$ that contains the two points $a^+_0$ and $a^-_0$; panel~(b) shows in $\Sigma$ all of $\hWu_0(\Gamma_1)$ (red curve), $\gamma^\pm_1$,  $\gamma^\pm_2$ and $a^\pm_0$.}
\end{figure}

\Fref{fig:Wu0} shows the primary intersection curve $\hWu_0(\Gamma_1)$, which is a single closed curve with two branches that connect $\gamma^\pm_1$ and $\gamma^\pm_2$; one of these two branches contains the points $a^+_0$ and $a^-_0$ and both branches cross the tangency locus $C$. Panel~(a) shows $\Gamma_1$ and the part of $W^u(\Gamma_1)$ that was computed to find the branch of $\hWu_0(\Gamma_1)$ containing $a^\pm_0$; starting from either of these points gives the same result. In contrast to \fref{fig:Ws}(a), which shows the projection onto $(c, c_t, n)$-space, the projection in \fref{fig:Wu0}(a) is now onto $(c, v, n)$-space. Hence, the plane $C$ appears as a line that divides the projection of $\Sigma$ into the regions $\Sigma^+$ and $\Sigma^-$; this projection illustrates that both $\Gamma_1$ and $\hWu_0(\Gamma_1)$ cross $C$. Panel~(b) shows all of $\hWu_0(\Gamma_1)$ in $\Sigma$; the second branch (not containing $a^\pm_0$) was found with continuation starting near $\gamma_1$. Note that $\hWu_0(\Gamma_1) \cap \Sigma^+$ maps to $\hWu_0(\Gamma_1) \cap \Sigma^-$ under the first-return map $f$; see also \cite{England2005,Lee2008}. Similar to what we found for $\hW_0^{s,\pm}(\Gamma_2)$, the intersection curve $\hWu_0(\Gamma_1)$ contains all backward images under $f^\pm$ of the two intersection points $a^\pm_0$, that is, $a^\pm_k \in \hWu_0(\Gamma_1)$ for all $k < 0$; note that all of these points (marked by stars)  are again extremely close to $\gamma^\pm_1$.

\begin{figure}[t!]
  \hspace*{1.5cm}
  \includegraphics{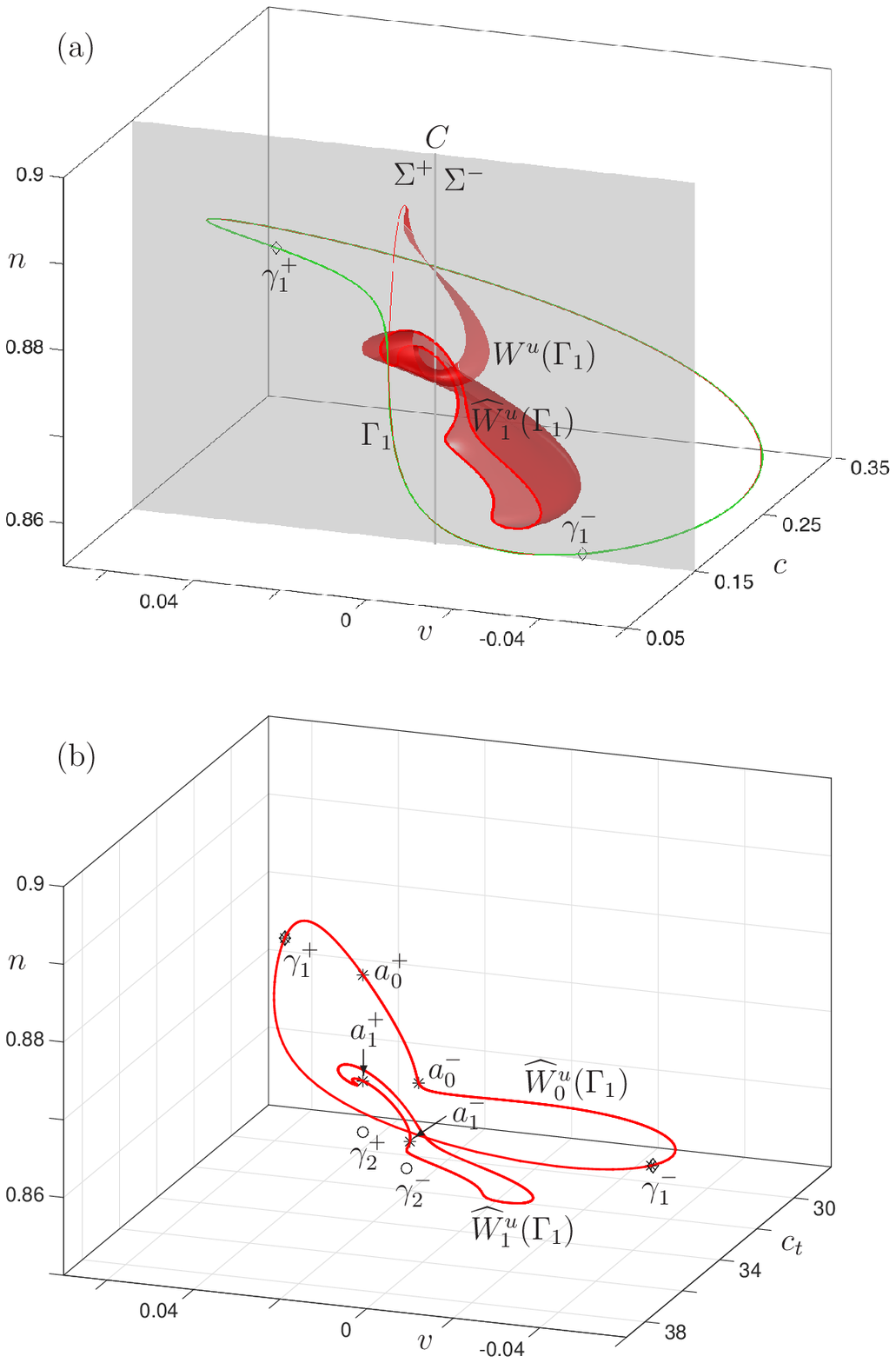} 
  \caption{\label{fig:Wu1}
    A part of $W^u(\Gamma_1)$ intersects $\Sigma$ (grey plane) again in the closed curve $\hWu_1(\Gamma_1)$. Panel~(a) shows, in projection onto $(c, v, n)$-space, the periodic orbit $\Gamma_1$ (green curve), the respective part of $W^{u}(\Gamma_1)$ (red surface) up to $\hWu_1(\Gamma_1)$ (red curve) in $\Sigma$; panel~(b) shows in $\Sigma$ the intersection sets $\hW^{u}(\Gamma_0)$ and $\hW^{u}(\Gamma_1)$ (red curves), $\gamma^\pm_1$,  $\gamma^\pm_2$, $a^\pm_0$ and $a^\pm_1$.}
\end{figure}

Since $\hWu_0(\Gamma_1)$ does not contain any of the points $a^\pm_k$ with $k > 0$, this means that $W^{u}(\Gamma_1)$ must intersect $\Sigma$ in other curves. In fact, most points on $\hWu_0(\Gamma_1)$ never return to $\Sigma$, but there is a small segment on $\hWu_0(\Gamma_1)$ that returns as the closed intersection curve $\hWu_1(\Gamma_1)$ and contains the points $a^\pm_1$; this small  segment on $\hWu_0(\Gamma_1)$ contains $a^+_0$ and maps under the first-return map $f$ to a similarly small segment containing $a^-_0$.

The curve $\hWu_1(\Gamma_1)$ can be found by continuation from either $a^+_0$ or $a^-_0$, and it is shown in \fref{fig:Wu1}. Panel~(a) shows, in projection onto $(c, v, n)$-space,  $\Gamma_1$ and the part of $W^{u}(\Gamma_1)$ that generates $\hWu_1(\Gamma_1)$. Panel~(b) shows $\hWu_1(\Gamma_1)$ and $\hWu_0(\Gamma_1)$ in $\Sigma$; also shown are $\gamma^\pm_1$,  $\gamma^\pm_2$, $a^\pm_0$ and $a^\pm_1$. Note from \fref{fig:Wu1}(a) that the closed curve $\hWu_1(\Gamma_1)$ has two intersections points with $C$, with $a^+_1$ and $a^-_1$ on either side, such that $\hWu_1(\Gamma_1) \cap \Sigma^+$ maps to $\hWu_1(\Gamma_1) \cap \Sigma^-$ under the first-return map $f$. Geometrically, the curve $\hWu_1(\Gamma_1)$ bounds a small ``bump'' of the surface $W^u(\Gamma_1)$ that crossed $\Sigma$ again. This type of intersection between a manifold and a global Poincar\'{e} section arises from a minimax quadratic tangency on $C$~\cite{Lee2008}. 

The intersection curve $\hWu_1(\Gamma_1)$ only contains $a^\pm_1$ and no other points in $A \cap \Sigma$. Hence, there must exist a further return of $W^u(\Gamma_2)$ to $\Sigma$, corresponding to an even smaller segment on $\hWu_0(\Gamma_1)$ that contains $a^+_0$, maps to a segment in $\Sigma^-$ that contains $a^-_0$, a segment in $\Sigma^+$ that contains $a^+_1$, one in $\Sigma^-$ that contains $a^-_1$, and finally, produces a new intersection curve $\hWu_2(\Gamma_1)$ that contains (at least) the point $a^+_2$. This curve $\hWu_2(\Gamma_1)$ can also be found by continuation starting from the orbit segment that ends on $a^+_2$.

\begin{figure}[t!]
  \hspace*{1.5cm}
  \includegraphics{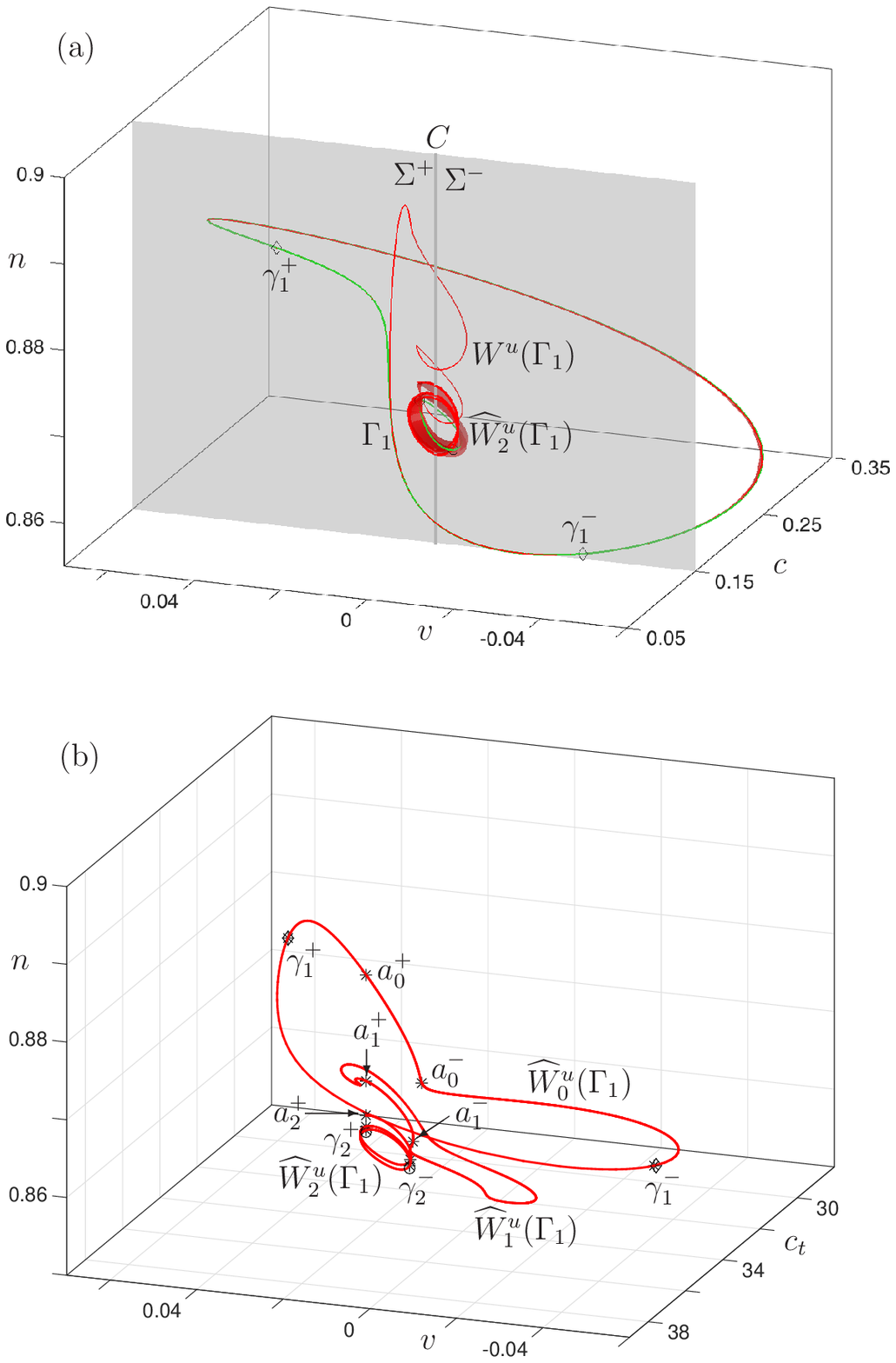} 
  \caption{\label{fig:Wu2}
    A part of $W^u(\Gamma_1)$ intersects $\Sigma$ (grey plane) a second time in a spiralling curve $\hWu_2(\Gamma_2)$ that contains the points $a^+_k$ for $k \geq 2$. Panel~(a) shows, in projection onto $(c, v, n)$-space, the periodic orbit $\Gamma_1$ (green curve) and the respective part of $W^{u}(\Gamma_1)$ (red surface) up to $\hWu_2(\Gamma_2)$  (red curve) in $\Sigma$; panel~(b) shows in $\Sigma$ the intersection sets $\hWu_0(\Gamma_0)$, $\hWu_1(\Gamma_1)$ and $\hWu_2(\Gamma_2)$ (red curves), $\gamma^\pm_1$,  $\gamma^\pm_2$, and $a^+_k$ for $k \geq 0$.}
\end{figure}

\Fref{fig:Wu2} shows the result of this continuation, in one direction, for the same projections as in \frefs{fig:Wu0} and~\ref{fig:Wu1}. Here, panel~(a) shows $\Gamma_1$ and the respective part of $W^u(\Gamma_1)$ that generates the computed part of $\hWu_2(\Gamma_1)$ starting at $a^+_2$, and panel~(b) shows $\hWu_2(\Gamma_1)$ with $\hWu_0(\Gamma_1)$ and $\hWu_1(\Gamma_1)$ in $\Sigma$. The curve $\hWu_2(\Gamma_1)$ spirals and crosses $C$ many times, which implies that it contains not only $a^\pm_2$ but also the points $a^\pm_k$ with $k > 2$. 

\begin{figure}[t!]
  \hspace*{1.6cm}
  \includegraphics{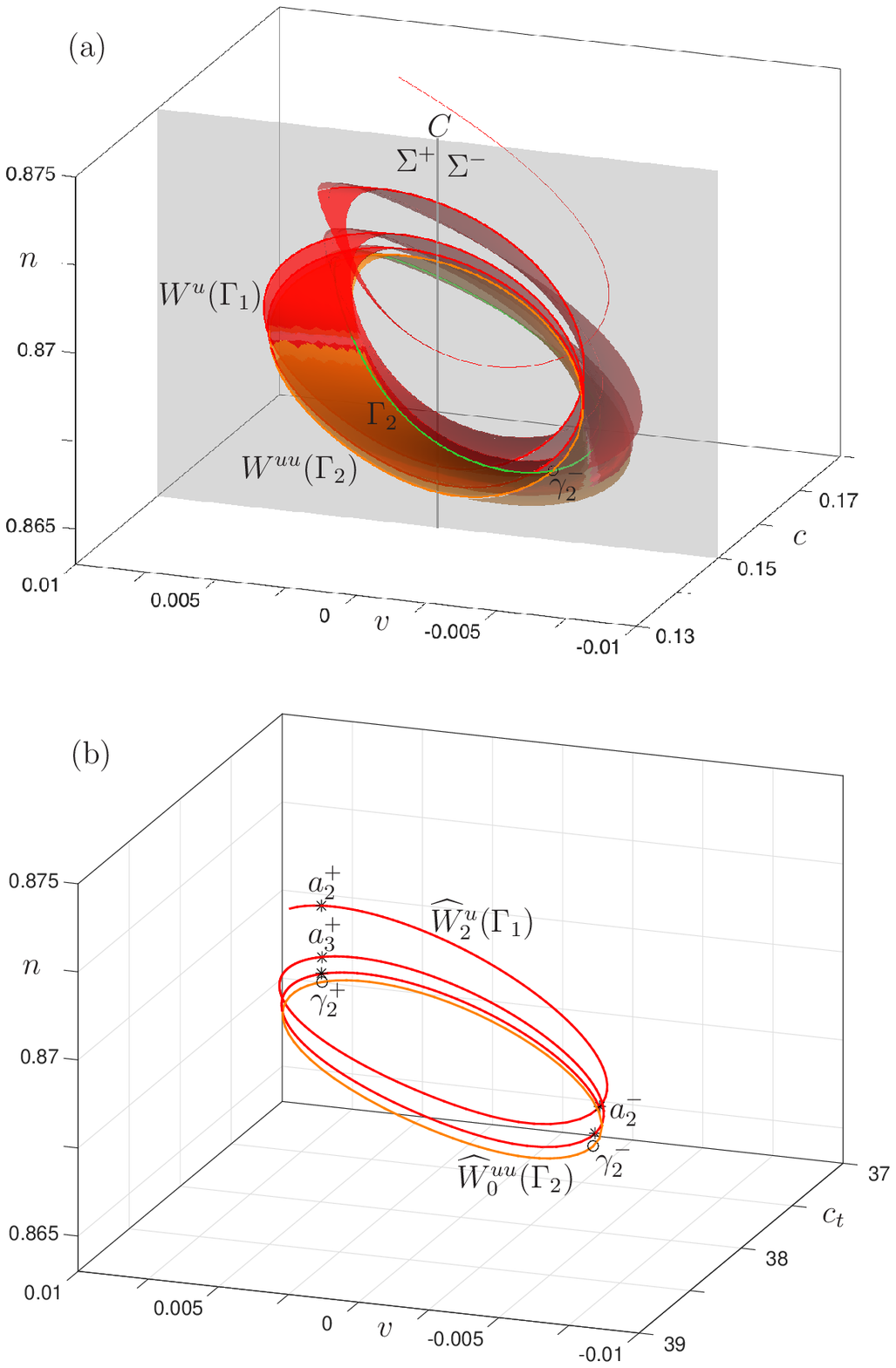}
  \caption{\label{fig:Wu2zoom}
    An enlargment of \fref{fig:Wu2} showing how $W^u(\Gamma_1)$ (red surface) spirals and accumulates onto the two-dimensional strong unstable manifold $W^{uu}(\Gamma_2)$ (orange surface). Panel~(a) shows a projection onto $(c, v, n)$-space with $\Sigma$ (grey plane); panel~(b) shows in $\Sigma$ the respective intersection sets $\hWu_2(\Gamma_0)$ (red curve),  $\hW_0^{uu}(\Gamma_2)$ (orange curve), $\gamma^\pm_2$, and $a^\pm_k$ for $k \geq 2$.}
\end{figure}

To see this better, \fref{fig:Wu2zoom} presents an enlargement. Its panel~(a) shows only $\hWu_2(\Gamma_1)$ and the local part of $W^u(\Gamma_1)$ that corresponds to this spiralling intersection curve. Furthermore, we also plot one side of the stong unstable manifold $W^{uu}(\Gamma_2)$ and its corresponding primary intersection curve $\hW_0^{uu}(\Gamma_2)$. The curves $\hWu_2(\Gamma_1)$ and $\hW_0^{uu}(\Gamma_2)$ are shown in $\Sigma$ in \fref{fig:Wu2zoom}(b). This panel clearly shows that $\hW_0^{uu}(\Gamma_2)$ is a single closed curve that crosses $C$ twice to connect the two points $\gamma^\pm_2$. Moreover, the curve $\hWu_2(\Gamma_1)$ accumulates on $\hW_0^{uu}(\Gamma_2)$ in a spiralling fashion. In particular, this means that $\hWu_2(\Gamma_1)$ contains all points $a^\pm_k$ for $k \geq 2$, which can be seen to converge to $\gamma^\pm_2$. When viewed locally near $\gamma^\pm_2$, this geometry is a result of the $\lambda$-lemma. On the other hand, the fact that all points  $a^\pm_k$ for $k \geq 2$ lie on $\hWu_2(\Gamma_1)$ is a global property related to how $W^u(\Gamma_1)$ intersects $\Sigma$.

\subsection{Interaction between $\hW^{s,\pm}(\Gamma_2)$ and $\hW^{u}(\Gamma_1)$}
\label{sec:wswu}

\begin{figure}[t!]
  \hspace*{1.5cm}
  \includegraphics{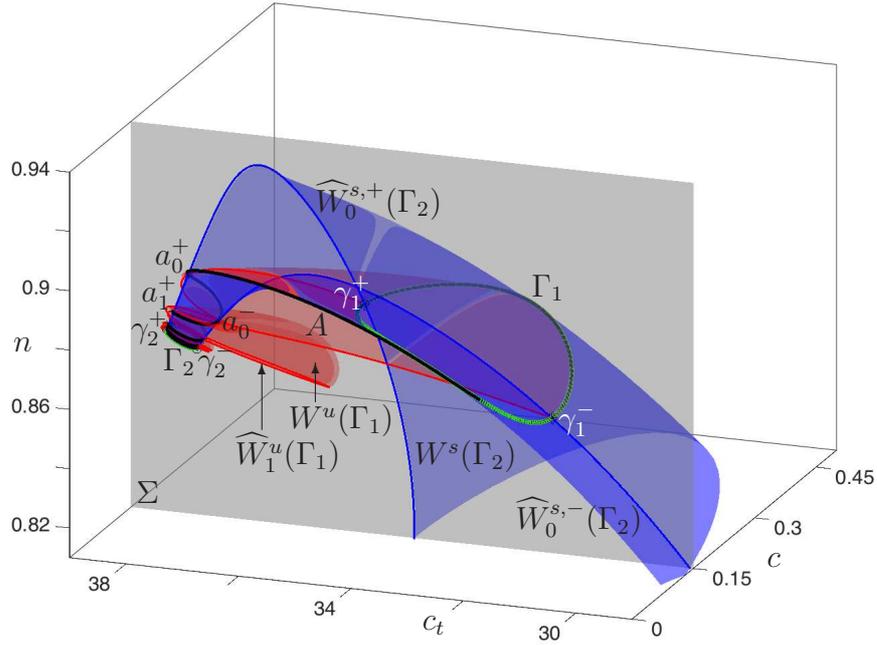}
  \caption{\label{fig:WsWu}
   An overall view in projection onto $(c, c_t, n)$-space of how $W^s(\Gamma_2)$ (blue surface) and $W^u(\Gamma_1)$ (red surface) intersect in the (non-transverse) connecting orbit $A$ (black curve), and how this generates the discrete intersection sets $a^\pm_k$ in the section $\Sigma$ (grey plane); compare with \fref{fig:Ws}(a).}
\end{figure}

Having accounted for all points in the sequence $\left( a_k^\pm \right)_{k \in \bbZ} = A \cap \Sigma$ as part of both $\hW^{s}(\Gamma_2)$ and $\hWu(\Gamma_1)$, we are now able to clarify how the codimension-one heterodimensional PtoP cycle arises from the intersection $W^{s}(\Gamma_2) \cap W^u(\Gamma_1)$. \Fref{fig:WsWu} illustrates, in projection onto $(c, c_t, n)$-space, how the two surfaces $W^s(\Gamma_2)$ and $W^u(\Gamma_1)$ intersect in the connecting orbit $A$ from $\Gamma_1$ to $\Gamma_2$. Also shown are the periodic orbits $W^{s}(\Gamma_2)$ and $W^u(\Gamma_1)$ and (the projection of) $\Sigma$; compare with \fref{fig:Ws}(a) and note the difference in projection in~\fref{fig:Wu0}(a). 

\begin{figure}[t!]
  \vspace*{2mm}
  \hspace*{1.5cm}
  \includegraphics{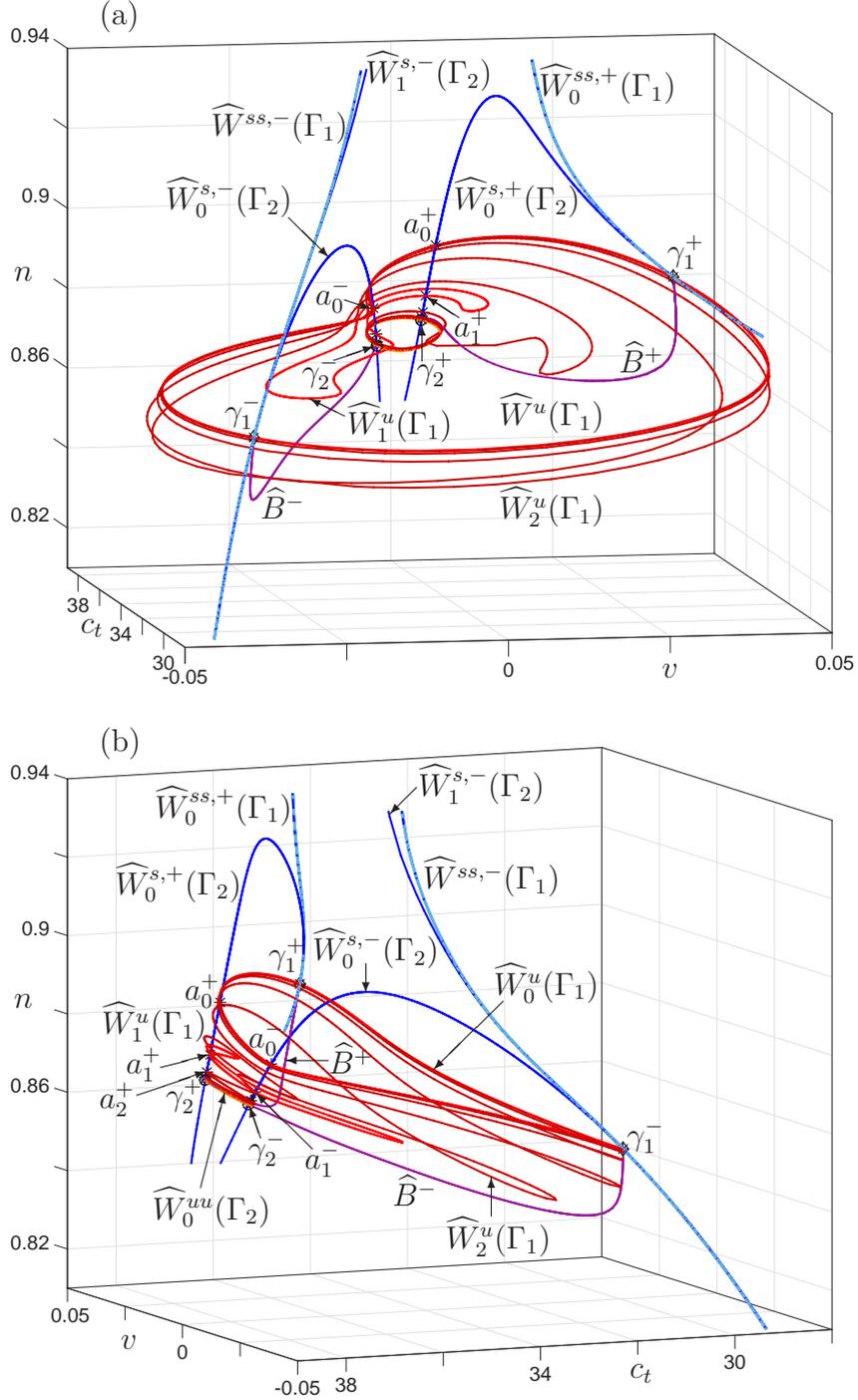}
  \caption{\label{fig:WsWuSelf}
    Two views of $\Sigma$ in panels (a) and (b) illustrate how $\hW^{s,+}(\Gamma_2)$ and $\hW^{s,-}(\Gamma_2)$ (blue curves) intersect $\hWu_0(\Gamma_1)$, $\hWu_1(\Gamma_1)$ and $\hWu_2(\Gamma_1)$ (red curves) in the points $a^+_k$ and $a^-_k$, respectively. Also shown are $\widehat{W}_0^{ss,\pm}(\Gamma_1)$, $\hW_0^{uu}(\Gamma_2)$ and the two intersection curves $\widehat{B}^\pm$ of the cylinder $B$; compare with \fref{fig:WsWu}.}
\end{figure}

\Fref{fig:WsWuSelf} assembles all relevant intersection sets in $\Sigma$ that are associated with the heterodimensional PtoP cycle. Shown are two different views featuring the curves in $\hW^{s,+}(\Gamma_2)$ and $\hW^{s,-}(\Gamma_2)$ as they intersect $\hWu_0(\Gamma_1)$, $\hWu_1(\Gamma_1)$ and $\hWu_2(\Gamma_1)$ in the heteroclinic points $a^\pm_k$ in the intersection set of $A$ with $\Sigma$. Also shown are the intersection curves $\widehat{W}_0^{ss,+}(\Gamma_1)$ and $\widehat{W}_0^{ss,-}(\Gamma_1)$ of the strong stable manifold near $\gamma_1^\pm$, and similarly, the intersection set $\hW_0^{uu}(\Gamma_2)$ of the strong unstable manifold near $\gamma_2^\pm$. Finally, \fref{fig:WsWuSelf} also shows the two intersection curves $\widehat{B}^\pm$ of the cylinder $B$ with $\Sigma$. 

Notice that \fref{fig:WsWuSelf} show a much larger part of the curve $\hWu_2(\Gamma_1)$, not only the part from \fref{fig:Wu2} that spirals onto $\hW_0^{uu}(\Gamma_2)$ and contains the points $a^+_k$ for $k \geq 2$. The curve $\hWu_2(\Gamma_1)$ keeps spiralling and has infinite arclength, even though it is confined to a bounded region. Note that $\hWu_2(\Gamma_1)$ also accumulates onto the closed curve $\hWu_0(\Gamma_1)$ that contains $\gamma_1^\pm$. Hence, $\hW^{u, \pm}(\Gamma_1)$ accumulates on itself, which implies the existence of homoclinic orbits to $\Gamma_1$. Notice further, that $\hWu_2(\Gamma_1)$ approaches $\gamma_1^\pm$ along the direction of $\widehat{B}^\pm$, that is, along the weak stable direction, as expected; see \fref{fig:Wu2}(b). 

\begin{figure}[t!]
  \hspace*{2.0cm}
  \includegraphics{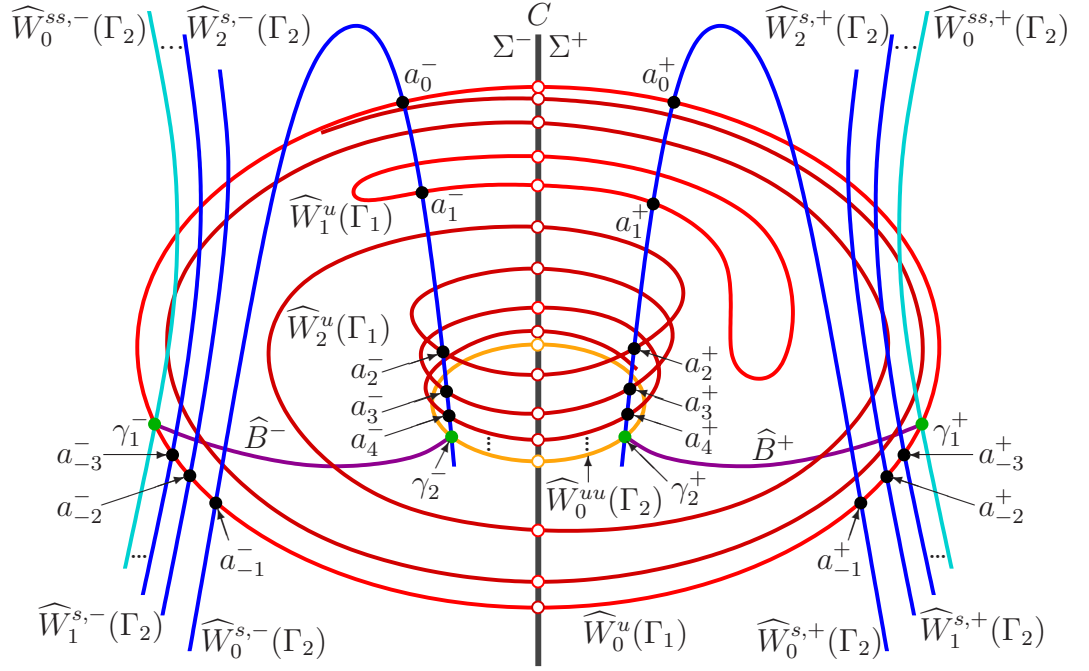}
  \caption{\label{fig:WsWuAllsketch}
    A sketch of the invariant objects in the section $\Sigma$ that give rise to the heterodimensional PtoP cycle; compare with \fref{fig:WsWuSelf}.}
\end{figure}

\Fref{fig:WsWuAllsketch} is a topological sketch of the objects shown in \fref{fig:WsWuSelf}. The tangency locus $C$ is indicated by a line that divides $\Sigma$ into the two parts $\Sigma^+$ and $\Sigma^-$. Illustrated is the accumulation of the curves $\hW^{s,\pm}(\Gamma_2)$ onto $\widehat{W}_0^{ss,\pm}(\Gamma_1)$, and the accumulation of $\hWu_2(\Gamma_2)$ onto $\hW_0^{uu}(\Gamma_2)$ and onto $\hWu_0(\Gamma_1)$. 

\Fref{fig:WsWuAllsketch} is designed to make the connection between the heterodimensional PtoP cycle in the four-dimensional continuous-time system~\eref{eq:4DAtri} and a heterodimensional cycle between two fixed points of a diffeomorphism in $\bbR^3$, as sketched in \fref{fig:hdcycle}. The first return $f$ on $\Sigma$ maps the curves in $\Sigma^+$ to those in $\Sigma^-$; while the second-return maps $f^\pm$ leave their respective domains $\Sigma^\pm$ invariant. We now restrict our attention to the action of $f^+$ on $\Sigma^+$; of course, the respective statements hold equally for $f^-$ on $\Sigma^-$. In \frefs{fig:WsWuSelf} and~\ref{fig:WsWuAllsketch}, the curves $\widehat{W}_0^{ss,+}(\Gamma_1)$ and $\widehat{B}^+$ are tangent at $\gamma_1^+$ to the strong stable and the weak stable direction, respectively, while the curve $\hWu_0(\Gamma_1)$ is tangent to the unstable direction. Similarly, at $\gamma_2^+$, the curves $\widehat{W}_0^{uu}(\Gamma_1)$ and $\widehat{B}^+$ are tangent to the strong unstable and the weak unstable direction, respectively, while the curve $\hW_0^{s,+}(\Gamma_1)$ is tangent to the unstable direction. This agrees exactly with the local situation sketched near $\gamma_1$ and $\gamma_2$ shown in \fref{fig:hdcycle}. Hence, the curve $\widehat{B}^+$ connects $\gamma_1^+$ with $\gamma_2^+$ in \frefs{fig:WsWuSelf} and~\ref{fig:WsWuAllsketch} exactly as the curve $\widehat{B}$ connects $\gamma_1$ with $\gamma_2$ in \fref{fig:hdcycle}. Moreover, the curves $\hW_0^{s,+}(\Gamma_2)$ and $\hWu_0(\Gamma_1)$ in \frefs{fig:WsWuSelf} and~\ref{fig:WsWuAllsketch}  intersect in the points $a^+_0$ and $a^+_{1}$ and also contains the points $a^+_k$ for $k \geq 2$, which is also as sketched in \fref{fig:hdcycle} for the curves $W^{s}(\gamma_2)$ and $W^u(\gamma_1)$. 

The difference is in the global nature of the interacting manifolds: although it is not shown in \fref{fig:hdcycle}, this sketch suggests that $W^{s}(\gamma_2)$ and $W^u(\gamma_1)$ are single curves that form a heteroclinic tangle by intersecting in the points $a_k$, where $W^{s}(\gamma_2)$ accumulates on $W^{ss}(\gamma_1)$ and $W^{u}(\gamma_1)$ accumulates on $W^{uu}(\gamma_2)$. This is the logical way of completing the picture for a diffeomorphism in $\bbR^3$. The diffeomorphism $f^+$ on $\Sigma^+$, on the other hand is not defined along the discontinuity locus $C$  and, moreover, certain parts of $W^{s}(\Gamma_2)$ and $W^{u}(\Gamma_1)$ do not intersect $\Sigma^+$; the latter results in the intersection sets $\hW_0^{s,+}(\Gamma_2) \cap \Sigma^+$ and $\hWu_0(\Gamma_1) \cap \Sigma^+$ in $\Sigma$ comprising of infinitely many curves. This reflects the fact that there exists no global section for system~\eref{eq:4DAtri}. \Frefs{fig:WsWuSelf} shows that both $W^{s}(\Gamma_2)$ and $W^{u}(\Gamma_1)$ accumulate on $W^{ss}(\Gamma_1)$ and $W^{uu}(\Gamma_2)$, respectively, which are locally cylinders. In the section $\Sigma$, this manifests itself differently near $\gamma_2^+$ and $\gamma_1^+$ because these local cylinders are intersected by $\Sigma$ differently: as two arcs in the case of $W^{ss}(\Gamma_1)$ and as a closed curve crossing $C$ in the case of $W^{uu}(\Gamma_2)$. Nevertheless, locally near the heterodimensional cycle of $f^+$ in $\Sigma^+$, formed by the point $\gamma_1^+$, the sequence $\left( a_k^+ \right)_{k \in \bbZ}$, the point $\gamma_2^+$, and the curve $\widehat{B}^+$, the manifold structure shown in \frefs{fig:WsWuSelf} and~\ref{fig:WsWuAllsketch} is, indeed, as sketched in \fref{fig:hdcycle}. 

Finally, \frefs{fig:WsWuSelf} and~\ref{fig:WsWuAllsketch} nicely illustrate why further heterodimensional PtoP cycles between $\Gamma_1$ and $\Gamma_2$ must exist for nearby parameter values. Namely, because of its accumulation onto $\hW_0^{u, \pm}(\Gamma_1)$, the curve of $\hW_2^{u, \pm}(\Gamma_1)$  comes arbitrarily close to $\hW^{s, \pm}(\Gamma_2)$. Hence, in any neighbourhood of the point $(J^\ast, s^\ast)$ there must be infinitely many further codimension-one PtoP connections, created by different intersections between the persisting one-dimensional sets $\hW^{s, \pm}(\Gamma_2)$ and $\hW^{u}(\Gamma_1)$. Since the cylinder $B$ is structurally stable near $(J^\ast, s^\ast)$, these will actually form further heterodimensional cycles, which involve more complicated excursions between $\Gamma_1$ and $\Gamma_2$. Moreover, the existence of homoclinic orbits to both $\Gamma_1$ and $\Gamma_2$ (see also \cite{ZKK2012}) implies the existence of many other periodic orbits and, hence, suggesting the existence of heterodimensional cycles involving different pairs of periodic orbits nearby in parameter space.

\section{Conclusions}
\label{sec:conclusions}

We showed how stable and unstable manifolds of two saddle periodic orbits with different indices intersect to form a heterodimensional PtoP cycle of a four-dimensional vector field. The PtoP cycle of the vector field intersects the transverse three-dimensional section $\Sigma$ twice on either side of a two-dimensional subspace $C$ along which the flow is tangent, which represents the typical situation for flows that are not periodically forced. Namely, the Poincar{\'e} map can be defined only locally (as a second return to $\Sigma$) near the respective intersection sets of the PtoP cycle on either side of $C$. We found that the structure of the manifolds agrees locally with that of two saddle fixed points in a three-dimensional diffeomorphism as considered in the literature on heterodimensional cycles. The global nature of the intersection set, on the other hand, is more complicated and involves intersection sets crossing the tangency locus $C$. Geometrically, the codimension-one part of the PtoP cycle in the four-dimensional phase space arises from the fact that: (i) the two-dimensional unstable manifold of the first periodic orbit accumulates on the two-dimensional strong unstable manifold of the second periodic orbit; and (ii) the two-dimensional stable manifold of the second periodic orbit accumulates on the two-dimensional strong stable manifold of the first periodic orbit. How these manifolds intersect the section $\Sigma$ near either of the two periodic orbits differs, and this explains the differences between the gobal properties of the respective intersection sets.

For these results we made extensive use of advanced numerical methods based on a two-point boundary value problem set-up to define and compute global (sub)manifolds of interest. This approach is versatile, efficient and, in particular, well suited to deal with large differences between expansion/contraction rates of the periodic orbits. It allowed us to compute and present the intersection sets of the respective manifolds in a three-dimensional section transverse to the two periodic orbits and to both heteroclinic connections that comprise the PtoP cycle. We also showed with three-dimensional projections how these intersection sets arise from the respective manifolds in the four-dimensional phase space.

Ongoing work concerns the structure of parameter space near the heterodimensional PtoP cycle. The codimension-one PtoP connection itself occurs along a curve in the bifurcation diagram in the $(J,s)$-plane, while the connecting cylinder is structurally stable. Our computations suggest that there should be further curves of PtoP cycles of the same two periodic orbits nearby, which correspond to different non-transverse intersections of their two-dimensional global manifolds. Finding and continuing these and other secondary hererodimensional PtoP connections and cycles is a considerable task. After all, theory predicts that the existence of heterodimensional PtoP can occur robustly, and the question is how this manifests itself in a concrete system. In particular, robustness of heterodimensional cycles requires that there are infinitely many periodic orbits of increasingly higher periods nearby. Hence, presenting an overall consistent picture in the parameter plane also requires the computation of other bifurcation curves, especially those of saddle-node bifurcations of periodic orbits and of their period-doubling bifurcations. Indeed, each of these curves may be an end point of a curve of PtoP connections. For example, we showed that the curve of the primary codimension-one PtoP connection ends on a curve of period-doubling bifurcation. While the primary codimension-one PtoP connection is then lost, the bifurcating period-doubled periodic orbit has the correct index to form a new heterodimensional cycle with the periodic orbit that does not bifurcate. Hence, this new PtoP connection could be seen as the continuation of the primary one. Moreover, we have evidence of the existence of a period-doubling cascade, suggesting that there may be an associated cascade of heterodimensional PtoP cycles to the successive period-doubled periodic orbits. A detailed investigation of the respective bifurcation sets in the $(J,s)$-plane near the primary PtoP connection will be presented elsewhere. 

From a more general perspective, our study shows that abstract concepts, introduced to advance the theory of dynamical systems, can be detected and investigated in concrete models, including in ODEs (generally the models of choice) rather than in diffeomorphisms; specifically, we considered the Atri model, which is from neuroscience and describes intracellular calcium dynamics. It is an interesting direction for future research to investigate how the complicated dynamics associated with a heterodimensional PtoP cycle manifests itself in experimental measurements. In case that there is evidence of a heterodimensional PtoP cycle playing an important role in organizing the system behaviour, the conclusion must be that the ODE model for the phenomenon studied requires a dimension of at least four.

\section*{Acknowledgements}

We thank Katsutoshi Shinohara for helpful discussions. BK and HMO were supported by Royal Society of New Zealand Marsden Fund grant 16-UOA-286.

%


\end{document}